\documentclass[onefignum,onetabnum,nonumbers]{siamart190516}
\pdfoutput=1

\usepackage{lmodern}
\usepackage{float}
\usepackage{amsfonts,amsmath,amssymb,graphicx,color}
\usepackage{amsmath}
\usepackage{verbatim}
\usepackage{color,graphicx}
\usepackage[algo2e,ruled,linesnumbered,vlined]{algorithm2e} 
\newtheorem{assumption}{Assumption}[section]  

\usepackage{cleveref}
\usepackage{subcaption}
\usepackage[normalem]{ulem}
\usepackage{enumitem}
\newcommand{\beq}{\begin{equation}}
\newcommand{\eeq}{\end{equation}}
\newcommand{\be}{\begin{equation}}
\newcommand{\ee}{\end{equation}}
\newcommand{\beqa}{\begin{eqnarray}}
\newcommand{\eeqa}{\end{eqnarray}}
\newcommand{\bea}{\begin{eqnarray}}
\newcommand{\eea}{\end{eqnarray}}
\newcommand{\beqas}{\begin{eqnarray*}}
\newcommand{\eeqas}{\end{eqnarray*}}

\newcommand{\eps}{\epsilon}

\newcommand{\bmat}{\left[ \begin{array}}
\newcommand{\emat}{\end{array} \right]}
\newcommand{\tpk}{\tilde{p}_k}
\newcommand{\tgk}{\tilde{g}_k}
\newcommand{\tpg}{\tilde{p}_k^T\tilde{g}_k}
\newcommand{\dgk}{\delta_g(x_k)}
\newcommand{\dpk}{\delta_p(x_k)}

\numberwithin{equation}{section}

\usepackage{dsfont}

\allowdisplaybreaks

\usepackage{lipsum}
\usepackage{amsfonts}
\usepackage{graphicx}
\usepackage{epstopdf}
\usepackage{algorithmic}
\ifpdf
  \DeclareGraphicsExtensions{.eps,.pdf,.png,.jpg}
\else
  \DeclareGraphicsExtensions{.eps}
\fi

\usepackage{enumitem}
\setlist[enumerate]{leftmargin=.5in}
\setlist[itemize]{leftmargin=.5in}

\newsiamremark{remark}{Remark}
\newsiamremark{hypothesis}{Hypothesis}
\crefname{hypothesis}{Hypothesis}{Hypotheses}
\newsiamthm{claim}{Claim}

\headers{Noise-Tolerant Optimization for Robust Design}{Y. Lou, S. Sun, J. Nocedal}

\title{
Design Guidelines for Noise-Tolerant Optimization with Applications in Robust Design\thanks{Submitted to the editors Jan 17th, 2024.
\funding{This work was supported by National Science Foundation grant  DMS-2011494, AFOSR grant FA95502110084, and  ONR grant N00014-21-1-2675.}}}

\author{
Yuchen Lou\thanks{Department of Industrial Engineering and Management Sciences, Northwestern University,  USA.} 
\and Shigeng Sun\thanks{Department of Engineering Sciences and Applied Mathematics,  Northwestern University,  USA.} 
\and Jorge Nocedal\footnotemark[2]
}

\usepackage{amsopn}

\begin{document}

\maketitle

\begin{abstract}
The development of nonlinear optimization algorithms capable of performing reliably in the presence of noise has garnered considerable attention lately. This paper advocates for strategies to create noise-tolerant nonlinear optimization algorithms by adapting classical deterministic methods. These adaptations follow certain design guidelines described here, which make use of estimates of the noise level in the problem. The application of our methodology is illustrated by the development of a line search gradient projection method, which is tested on an engineering design problem. It is shown that a new self-calibrated line search and noise-aware finite-difference techniques are effective even in the high noise regime. Numerical experiments investigate the resiliency of key algorithmic components. A convergence analysis of the line search gradient projection method establishes convergence to a neighborhood {of stationarity}.
\end{abstract}



\section{Introduction}
\label{sec:intro}
Over the past 50 years, significant progress has been made in the development of robust and efficient methods for nonlinear optimization. 
These methods have been adopted in a wide range of applications, and in the case of constrained optimization, can be quite complex. Recently, there has been a growing interest in tackling nonlinear problems where the function and/or gradient evaluations are subject to noise or errors
\cite{toint2021TRqEDAN,berahas2021stochastic,cartis2021TRqDA,katya2018storm,curtis2019stochTR,toint2021ARqpEDA2,more2012estimating,ng2014multifidelity,xie2020analysis}. 
 This raises the question of whether existing optimization methods require substantial redesign to ensure robustness in the presence of noise, or if certain modifications are sufficient to tackle such challenges.

This paper argues that one can develop effective methods for a broad range of noisy optimization problems by retaining the fundamental properties of deterministic methods while incorporating certain modifications based on the design guidelines outlined herein.
These guidelines stem from the observation that, in the presence of noise, only  few operations can lead to numerical difficulties in optimization methods. These operations include:
\begin{enumerate}
\item Comparisons of noisy function values, as required e.g., in line search and trust region techniques. 
\item  Computation of differences of noisy function values, as required in finite-difference approximations to a gradient. 
\item  Computation of differences of noisy gradients, a basic ingredient in quasi-Newton updating.
\end{enumerate}

{
Robust methods can be designed by ensuring that these operations are conducted reliably, preventing the algorithm from making harmful decisions. In this paper, we explore stabilization procedures that utilize an upper bound or a standard deviation of the noise (referred to as the \emph{noise level}), and illustrate their performance in solving a design optimization problem.
Examples of strategies proposed in the literature  to safeguard the three fragile operations mentioned above are as follows.  
{\textit{Soft comparisons:}} when assessing whether a step is acceptable by comparing noisy function values, the classical sufficient decrease condition can be relaxed in proportion to the noise level \cite{berahas2019derivative,berahas2019global,sun2022trust};
{\textit{Robust difference intervals:}} in computing a finite difference gradient approximation, the distance between evaluation points for noisy functions should be proportional to the square root of the noise level divided by the norm of the true Hessian \cite{more2012estimating,shi2022adaptive};
{\textit{Controlled gradient differences:}} quasi-Newton methods can achieve robustness by ensuring that points used for computing gradient differences (normally consecutive iterates) are adequately spaced in relation to the noise level in the problem \cite{shi2022noise,xie2020analysis}. 

We do not argue  that \emph{the only} way to design nonlinear optimization methods for noisy problems is to adapt existing deterministic methods. 
We will see that in scenarios with highly noisy gradients, deviating from traditional approaches can be beneficial. Specifically, utilizing techniques like diminishing steplengths \cite{bertsekas2011incremental,nemirovski2009robust,RobMon51} can help counteract the adverse impacts of errors or noise, offering a viable alternative to line searches or trust region techniques.
Nevertheless,  the sophistication of some of the established methods and software for deterministic optimization makes it alluring to build upon their foundations as much as possible because of the important algorithmic ideas they embody.
 For example, in cases where a good estimate of the optimal active set is available, it is sensible to employ an active set method like sequential quadratic programming, as it can effectively utilize this estimate \cite{robinson1974perturbed}. Similarly, primal-dual interior point methods have demonstrated remarkable efficacy in handling large-scale problems with network structure \cite{granville1994optimal}. Maintaining these capabilities even amidst noise is highly desirable.

In this paper, we study the performance of an algorithm that follows the design principles mentioned above and apply it to a design optimization problem in which the noise level can be adjusted. 
In this problem, the goal is to optimize the shape of an acoustic horn to achieve optimal efficiency, assuming that there is uncertainty in some of the physical properties of the system \cite{ng2014multifidelity}. 
This leads to a nonconvex bound constrained optimization problem, for which we design a noise-tolerant gradient projection method with a new \emph{self-calibrated line search} that adapts classical line search parameters to cope with  noise. 
Our case study provides ample flexibility for assessing the efficacy of various optimization methods as noise increases from mild  to extremely high, a regime where the stochastic gradient descent (SGD) method \cite{RobMon51} has shown to be particularly effective.

\subsection{Contributions of the Paper}
The recent literature on noisy nonlinear optimization typically reports numerical tests using either synthetic noise or simple machine learning models, leaving the question of their effectiveness in realistic applications open. In this paper, we focus on the sources of noise and errors that arise in certain practical problems, identify three critical operations prone to failure, and  discuss the importance of the noise level in designing noise-tolerant algorithms. Based on a case study in optimal design, we conduct systematic tests to verify the robustness of two key components of our gradient projection method, namely the line search and the finite difference gradient approximation, as the noise level in the problem increases.

Building upon these findings, we introduce a new \emph{self-calibrated line search} technique, effective even in environments with high levels of noise. This technique narrows the gap between traditional algorithms and the fixed step length projected SGD method. Additionally, we provide a convergence analysis for the line search gradient projection algorithm used in our case study, under the assumption that the noise in the function is bounded---a realistic assumption in this context.

 \subsection{Organization of the Paper}

This paper is structured into seven sections. In the following section, we explore the concept of \emph{noise level} and its estimation. Section \S\ref{sec:case} introduces the optimal design problem central to our study. In Section \S\ref{sec:gradp}, we detail a gradient projection method rooted in robust design principles. Section \S\ref{sec:numerical} presents the results of our numerical tests, while Section \S\ref{sec:convergence} offers a global convergence analysis of the gradient projection method with a line search. The paper concludes with final remarks in Section \S\ref{sec:finalr}.

\section{Noise and Errors}   \label{sec:ne}

Let {$f: \mathbb{R}^n \rightarrow \mathbb{R}$} be a smooth function and  $\tilde f$   its noisy or inexact counterpart.  Polyak \cite{polyak1987introduction} proposed two broad categories of noise and errors:
\begin{equation}
 \tilde f(x)  =  f(x) + \Delta(x) \quad\mbox{stochastic noise},  \label{probs}
  \end{equation}
 \begin{equation}  \tilde f(x)  =  f(x) + \delta(x) \quad\mbox{deterministic error}.  \label{probd} 
\end{equation}
The first case arises e.g. from Monte Carlo simulation, {and thus} $\Delta(x)\sim D_x$ is a random variable {following a distribution $D_x$ that may be parameterized by $x$}. 
The second case concerns computational error, broadly speaking, where repeated evaluations of $\tilde f(x)$ for a given $x$ give the same result.

Following Mor\'e and Wild \cite{more2011estimating,more2012estimating}, we use the term \emph{noise level} of a function. 
For the case of stochastic noise, we define the noise level of  $\tilde f$ at a point $x$ {as the standard deviation of $\tilde f(x)$, which we denote $\sigma_f(x)$.} In practice, we compute   an estimate $\epsilon_f(x)$:
 \begin{equation}   
 \label{eq:noise_level}
 \epsilon_f(x)\approx\sigma_f(x):=
 \sqrt{\mathbb{V}(\tilde f(x))},
 \end{equation} 
{where $\mathbb{V}(\cdot)$ denotes the variance of a random variable.
We argue that~\eqref{eq:noise_level} is adequate for our purposes,  because even though $\Delta(x)$ may contain bias, the bias term will cancel out in the three operations listed above, which involve comparisons or differences between function values.
We will return to this in Section~\ref{sec:choosing_epsA}.
}

    There are situations where deterministic error  \eqref{probd} can be described in a useful manner using a {stochastic} model, {so that  $\delta(x)$ can be viewed as a realization of a random variable}.
    In this case, we say that the function exhibits  \emph{computational noise},  and we will denote the resultant random variable as $\Delta(x)$,
    as in the case of stochastic noise. Following  Mor\'e and Wild \cite{more2011estimating,more2012estimating}, we 
     define the noise level $\sigma_f(x)$ as the standard deviation of $\Delta(x)$, with $\epsilon_f(x)$ serving as an approximate measure.
    For example, roundoff error is deterministic but can be modeled (albeit imperfectly) using a random variable  drawn from a uniform distribution over the interval $[-|f(x)|\epsilon_M,\, |f(x)|\epsilon_M]$, where $\epsilon_M$ is unit roundoff. 
    More examples of computational noise can be found in  \cite{more2011estimating} and in \S\ref{sec:compn} of this paper.

    In summary, stochastic and computational noise can be analyzed using a uniform approach by studying the properties of $\Delta(x)$.

In the case of {canonical} deterministic error, we can employ an estimate of the maximum error:
 \begin{equation}  \label{ebound}
     \epsilon_b \approx \sup |\delta(x)|, \quad x \in {\cal R}, 
\end{equation}
     where ${\cal R}$ is the region of interest {where one would expect the optimization process to take place}. 
     {In this paper, we assume the existence of a bound or a standard deviation of the noise.}

\subsection{Noise Level Estimation}
\label{sec:noise} 
Knowledge of the noise level in the function is a key component in the algorithms described in this paper. As a result, we now discuss some practical procedures for estimating the noise level.


\smallskip\noindent\textit{ Local Pointwise Estimate} ${ \epsilon_f(x).}$  
Given $m$ i.i.d. samples  $\{\tilde f_1(x), \tilde f_2(x),\ldots,\tilde f_m(x)\}$, we can  define the pointwise noise level, in the case of stochastic noise, as
\be
\label{eq:noise_level_pointwise}
\epsilon_f(x):=\sqrt{\frac{1}{m-1}\sum_{j=1}^m \left(\tilde f_j(x)-\overline{\tilde f(x)}\right)^2}, \quad \mbox{where} \ \  \overline{\tilde f(x)}:=\frac{1}{m}\sum_{j=1}\tilde f_j(x).
\ee
From classic statistics, we know that $\epsilon_f(x)$ is {a consistent estimator} of $\sigma_f(x)= [ \mathbb{V}(\tilde f(x))]^{1/2}$.

We observe that formula \eqref{eq:noise_level_pointwise} is not suitable in the context of computational noise. Since this type of noise is deterministic, the formula would erroneously suggest a noise level of zero.
One can, however, use the \texttt{ECNoise} algorithm  \cite{more2011estimating}, which was specifically designed for computational noise. It samples points along a randomly chosen line and employs Hamming differences \cite{hamming2012introduction} to yield an estimate $\epsilon_f(x)$.

\smallskip
\noindent\textit{Global Estimate $\epsilon_f$.}
Estimating $\epsilon_f(x_k)$ at every iteration is expensive and often unnecessary in practice. Whenever possible, it is desirable to employ a universal estimate $\epsilon_f$  for all $x$ in the region of interest. A global measure of noise over the region of interest ${\cal R}$, {which is assumed to be bounded,} can be defined as
\be
\overline \sigma_f = \frac{1}{|{\cal R}|}\int_{{\cal R}}\sigma_f(x)dx,
\ee
and can be estimated as
\be
\label{eq:epsf_avg}
{\epsilon}_f := \frac{1}{M}\sum_{i=1}^M\epsilon_f(x_i) \approx \overline \sigma_f
\ee
where $\{x_1,\cdots,x_M\}$ are  randomly sampled from ${\cal R}$, and $\epsilon_f(x_i)$ {is either given by \eqref{eq:noise_level_pointwise} or is the output of \texttt{ECNoise}}. 

In some cases, e.g.  Figure~\ref{fig:std100} in the next section, $\sigma_f(x)$ remains relatively constant across ${\cal R}$, allowing us to use a few (ideally only one representative) sample point $x_i$ to define $\epsilon_f$.


There are other more powerful estimators in the statistics literature  but they are typically more expensive. Considering the iterative aspect of optimization algorithms, the simpler constant estimators $\epsilon_f$ defined above are often adequate for practical purposes, as illustrated in \S\ref{sec:numerical}.
\section{Case Study: An Acoustic Design Problem} \label{sec:case}
To guide our discussion on the design of robust optimization methods and illustrate the concept of noise level, we begin by presenting a case study involving optimal design under uncertainty. In this problem, the uncertainty of some system parameters and the use of sampling techniques lead to noise in the objective function. While the uncertainty in the parameters is well-defined, predicting its  propagation into the objective function becomes challenging owing to the nonlinear nature of the simulation. Nonetheless, we will see that estimating the noise level in the function is feasible, enabling us to effectively utilize a range of approaches to solve the optimization problem.

\subsection{Statement of the Problem}

We consider the 2-D acoustic design problem under uncertainty studied by Ng and Willcox \cite{ng2014multifidelity}.
An incoming wave enters a horn through its inlet and exits the outlet into the exterior domain with an absorbing boundary; see Figure~\ref{fig:horn}. 
The goal is to find the shape of the horn so as to optimize its efficiency.

\begin{figure}[!htbp]
\begin{center}
\includegraphics[width=0.4\textwidth]{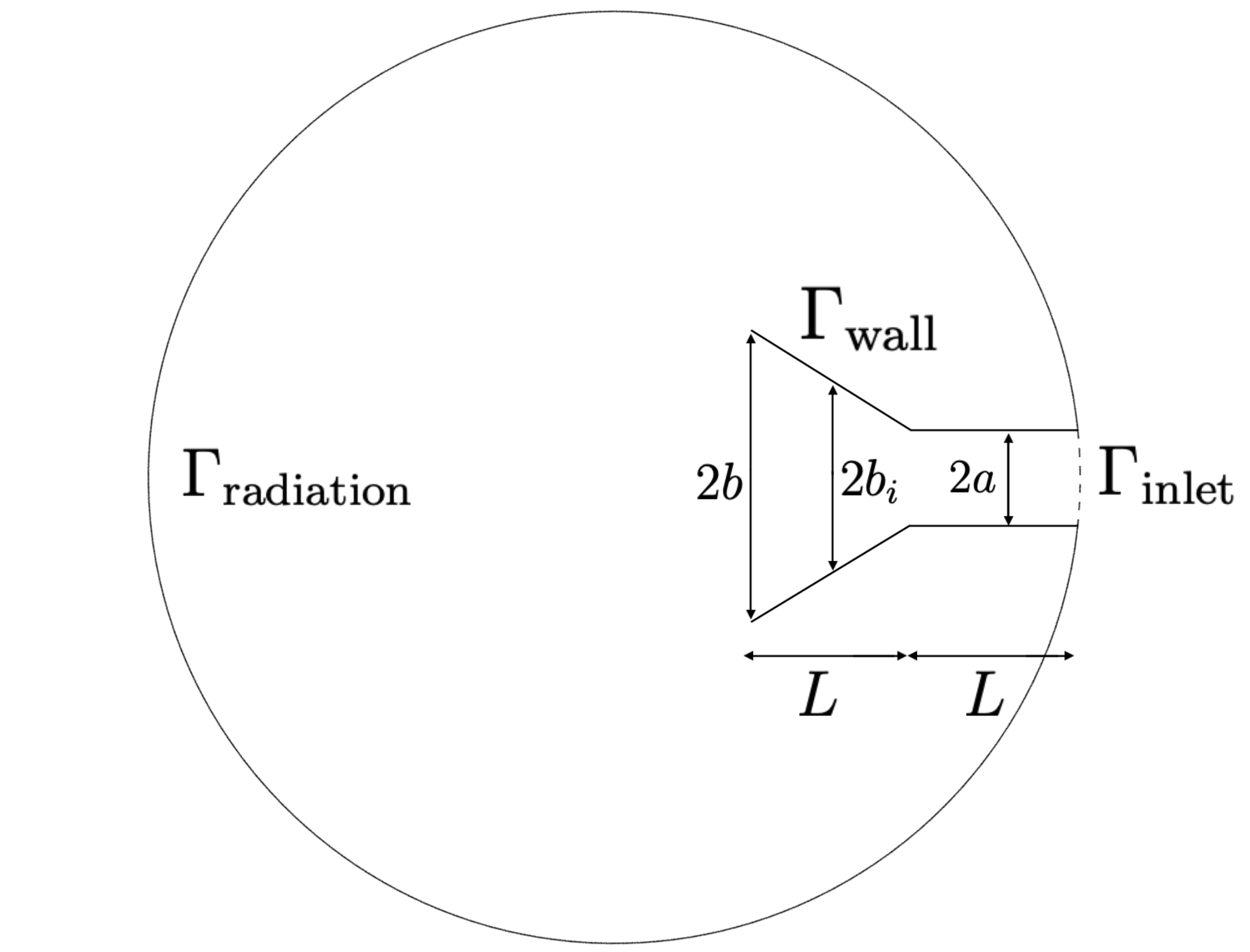}
\caption{Schematic plot for the design of horn}
\label{fig:horn}
\end{center}
\end{figure}

The propagation of the acoustic wave is modeled by the {dimensionless} complex Helmholtz equation 
\be 
\label{eq:pde}
\nabla^2 u+ \hat k^2u=0,
\ee where $u$ represents velocity and $\hat k$ is the wave number. The design variables $b=(b_1,b_2,\cdots,b_6)$ in $\mathbb{R}^6$ define the flare half-widths. We impose  bounds on the design variables,  $b_L \leq b \leq b_U$, and  assume that the dimensions, $a, L$,  depicted in Figure~\ref{fig:horn} are given. {The homogenous Neumann boundary condition is prescribed at the horn wall $\Gamma_{\text{wall}}$, and the first-order (Sommerfeld) radiation boundary condition is applied at $\Gamma_{\text{radiation}}$ and a source at the inlet $\Gamma_{\text{inlet}}$. For more details, see \cite[Section 2.2.2]{eftang2012two} and \cite[(4)]{udawalpola2008optimization}.} The PDE is solved using a finite element method.

The model contains uncertainties. The impedances $z_l$ and $z_u$ of the lower and upper horn walls are not known, but are assumed to follow a Gaussian distribution,  $N(50,3)$. Similarly, the wave number $\hat k$ is   assumed to follow a uniform distribution Unif$(1.3,1.5)$. We characterize uncertainty by the random variable $\omega$, so that $\hat k(\omega)\sim\text{Unif}(1.3,1.5)$, $z_l(\omega)\sim N(50,3)$, and $z_u(\omega)\sim N(50,3)$.

For a particular realization $\xi_i$ of the random variable $\omega$, the efficiency $s$ of the horn is characterized by the flux at the inlet, as follows:
\begin{equation}
\label{eq:efficiency}
    s(b, \xi_i) =\left|\int_{\Gamma_{\text{inlet}}}u (b, \xi_i, t)  dt-1\right|.
\end{equation}
Ng and Wilcox  employ various statistics of $s(b,\omega)$ to estimate overall efficiency and to achieve a robust design. We focus here on the following formulation
\begin{equation}
\label{prob:horn}
    \min_{b_L\leq b \leq b_U} f(b)=\mathbb{E}[s(b,\omega)]+3\sqrt{\mathbb{V}[s(b,\omega)]}.
\end{equation}
Although one may argue in favor of other robust formulations, the precise choice of the objective is not important in the discussion that follows. Note that problem \eqref{prob:horn} is a bound constrained stochastic optimization problem. 

\subsection{Approximating the Objective Function}
Closed form representations of the expectation and variance terms in \eqref{prob:horn} are unknown and must be estimated by sampling. 
At every iteration $k$  of the optimization algorithm, we compute the stochastic approximation:
 \begin{equation} \label{precisely}
      \tilde f(b_k)= \bar{s}_k(b_k,\Xi_k)+3\sqrt{S_{k}(b_k,\Xi_k)^2} ,
      \end{equation}
 where $\Xi_k=\{\xi_1,\xi_2,\cdots,\xi_N\}$  is a batch of i.i.d. samples of the random variable $\omega$. Here,   $\bar{s}_k(b_k,\Xi_k)$ is the sample mean of $s(b_k,\xi_i)$ with respect to the batch $\Xi_k$, i.e., 
\begin{equation}  \label{smalls}
    \bar{s}_k(b_k,\Xi_k)=\frac{1}{N}\sum_{\xi_i\in \Xi_k} s(b_k,\xi_i),
\end{equation}
and  $S_{k}(b_k,\Xi_k)^2$ is the sample variance of $s(b_k,\xi_i)$ in $\Xi_k$, {i.e.},
\begin{equation}  \label{larges}
    S_{k}(b_k,\Xi_k)^2=\frac{\sum_{\xi_i\in \Xi_k}\left(s(b_k,\xi_i)-\bar{s}_k(b_k,\Xi_k) \right)^2}{N-1}.
\end{equation}
For simplicity, we assume the batch size  $|\Xi_k|=N$ is constant across all optimization iterations.

 The evaluation of $\tilde f$ is expensive because, for each of the $N$ realizations of $\omega$, the acoustic efficiency $s$ given in \eqref{eq:efficiency} requires the solution of a differential equation using a finite element method that involves the solution of a linear system of order $\mathcal{O}$({30,000}). 
 (Ng and Willcox \cite{ng2014multifidelity}  employ a multifidelity approach to improve the efficiency of the sampling mechanism, but we will not consider it as it is not central to this investigation.)
\subsection{Illustration}
To visualize the behavior of the noisy function \eqref{precisely}, we plot it in Figure~\ref{fig:obj100} over a two-dimensional slice of $\mathbb{R}^6$ defined by varying two variables: $b_3$, $b_4$. 
The noise displays a discernible pattern rather than being highly erratic. As a result, the optimization problem is tractable notwithstanding the inherent nonlinearity of the simulation.

\begin{figure}[htpt!]
    \centering
    \includegraphics[width = 0.32\textwidth]{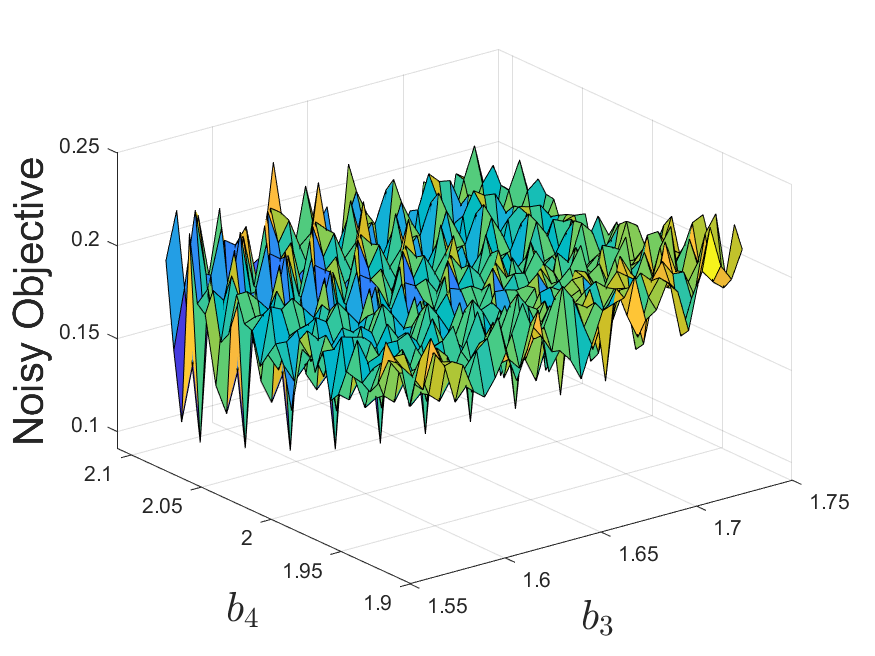}
    \includegraphics[width = 0.32\textwidth]{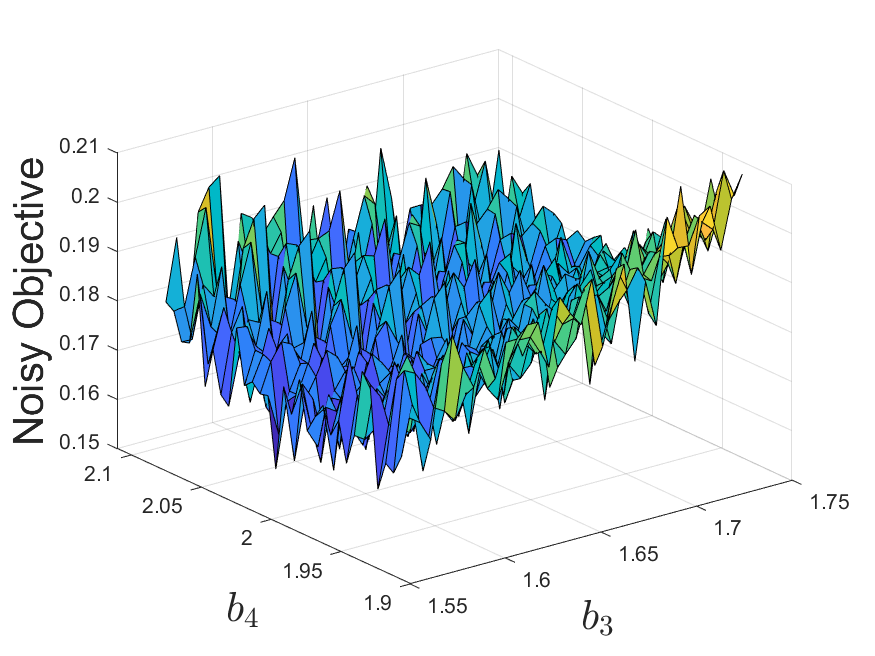}
    \includegraphics[width = 0.32\textwidth]{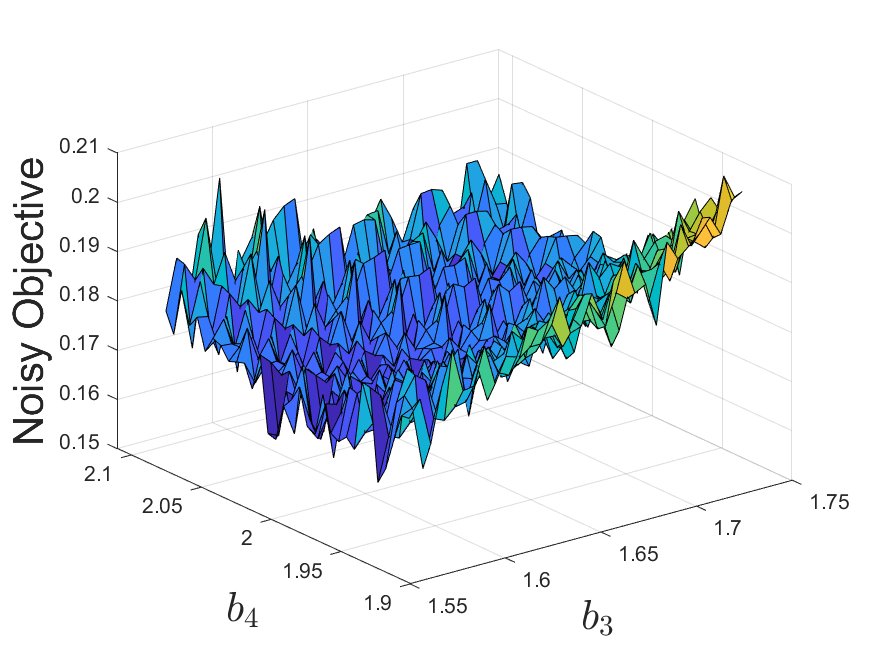}
    \caption{Noisy Function. The vertical axis plots the noisy objective \eqref{precisely} with different numbers of sample points: $N=10$ (left), $N=50$ (middle), and $N=100$ (right). The horizontal axes represent values of two of the design variables, {$b_3$ and $b_4$}. Different realizations of the random variable {$\omega$} were employed for each evaluation of $\tilde f$ in the region of interest.}
    \label{fig:obj100}
\end{figure}
The noise level $\sigma_f(b)$ in this function is defined as the standard deviation of $\tilde{f}(b)$ (see \eqref{eq:noise_level}), since  the problem exhibits stochastic noise.      
In Figure~\ref{fig:std100}, we plot an estimate $\epsilon_f(b)$ of $\sigma_f(b)$  (defined in~\eqref{eq:noise_level_pointwise} with $m=50$) as we vary  the variables $b_3, b_4$. While $\epsilon_f(b)$ does vary among different values of $b$, its fluctuations are not substantial.  
Thus, a single estimate might suffice for the optimization, as discussed in sections below.

\begin{figure}[!htbp]
    \centering
    \includegraphics[width = 0.32\textwidth]{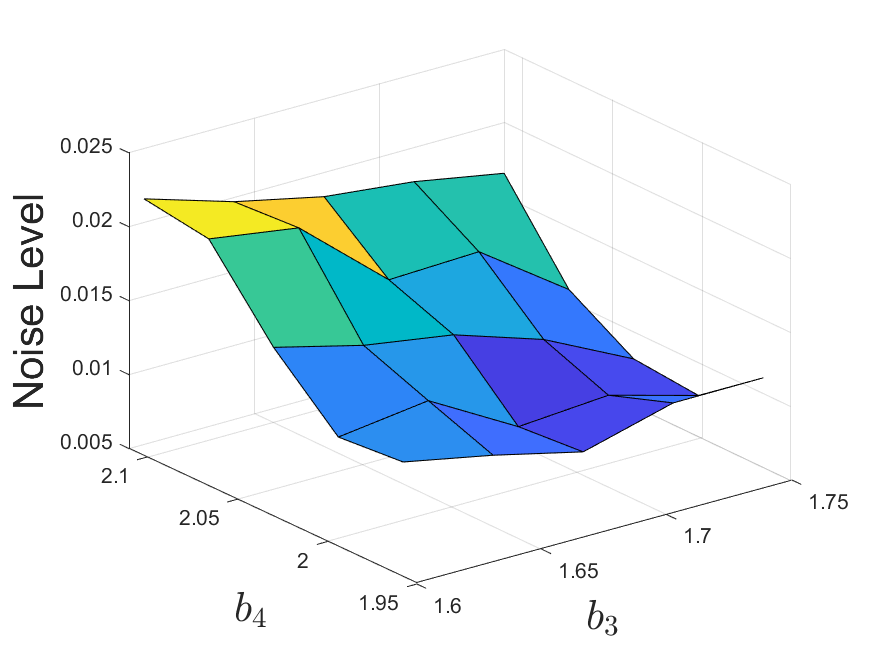}
    \includegraphics[width = 0.32\textwidth]{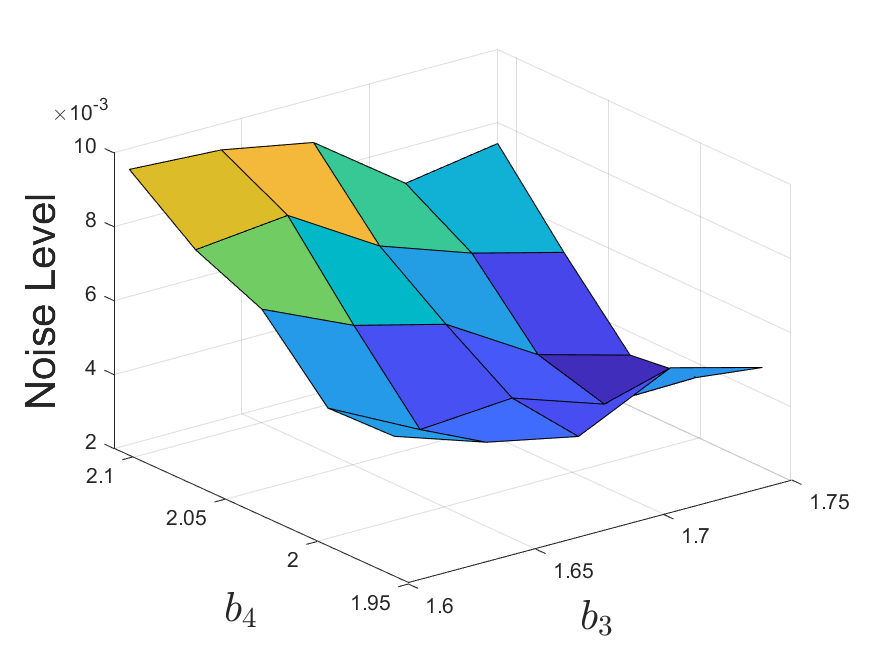}
    \includegraphics[width = 0.32\textwidth]{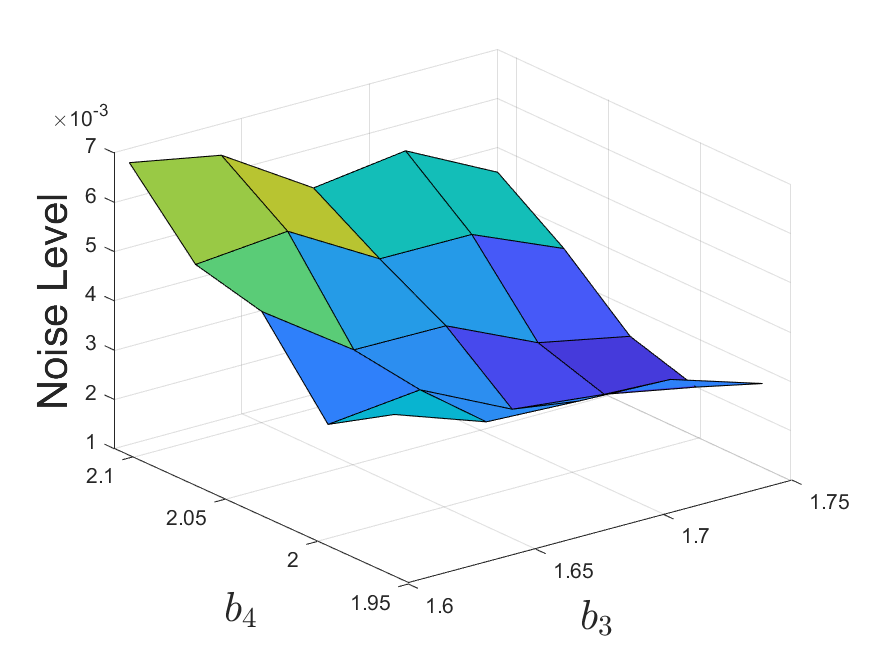}
    \caption{Noise Level. The vertical axis plots the estimated noise level $\epsilon_f(b)$ of the objective \eqref{precisely} with different numbers of sample points: $N=10$ (left), $N=50$ (middle), and $N=100$ (right). The horizontal axes represent values of two of the design variables $b_3$ and $b_4$. Each value $\epsilon_f(b)$ is computed as defined in~\eqref{eq:noise_level_pointwise} with $m=50$.}
    \label{fig:std100}
\end{figure}

\subsection{A Variant Illustrating Computational Noise} 
\label{sec:compn}

The acoustic horn design problem can also be used to illustrate computational noise. As already mentioned, the finite element solution of the Helmholtz equation \eqref{eq:pde} requires solving of a non-symmetric linear system of equations. In the examples given in the previous sections this was done using a direct linear solver. However,  practical applications often benefit from approximating solutions with an iterative method. 
In our next experiment, we utilize the {\sc gmres} method, with tolerance of $10^{-6}$, to solve the linear system.

 In order to isolate the effect of computational noise, we generate and fix a particular realization of $\Xi_k\equiv \Xi$, of size $N=10$, in the evaluation of the objective function \eqref{precisely}. We plot the generated objective function in Figure~\ref{fig:comp_noise}. {For comparison}, we also plot the function using a direct linear solver.
 In this context, computational noise is notably smaller than the stochastic noise previously illustrated. Although increasing the linear solver's tolerance can amplify the noise level, for brevity, our experiments will concentrate solely on stochastic noise.

\begin{figure}[!htbp]
    \centering
    \includegraphics[width = 0.40\textwidth]{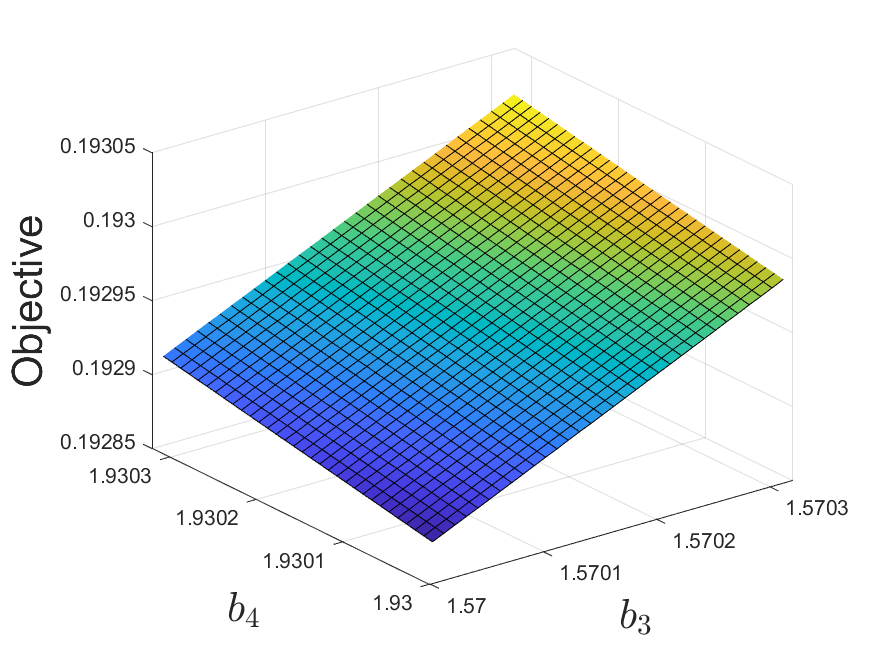}
    \includegraphics[width = 0.40\textwidth]{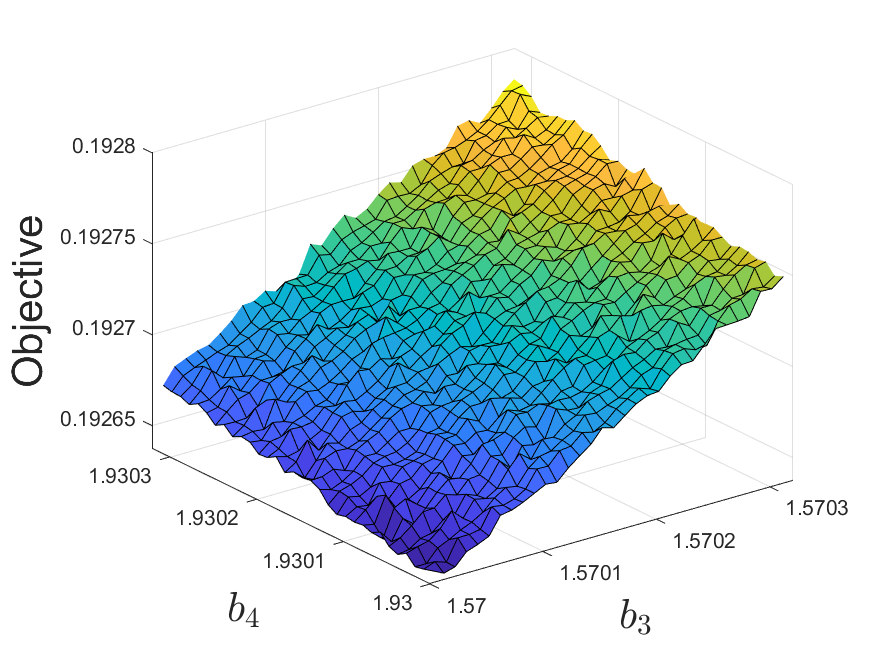}
    \caption{Computational Noise. The vertical axis plots a deterministic variant of the function~\eqref{precisely} in which the samples have been fixed. The linear system within the PDE scheme is solved using a direct method (left)  and using the iterative method {\sc gmres} with tolerance $10^{-6}$   (right).} 
    \label{fig:comp_noise}
\end{figure}


\section{Line Search Gradient Projection Methods}
\label{sec:gradp}
In this section, we consider algorithms for the solution of noisy bound constrained  optimization problems, such as the acoustic design problem \eqref{prob:horn}. Our starting point is a classical gradient projection method with a backtracking line search, designed to be stable with respect to the critical operations discussed in the introduction.  

Let $g(x):=\nabla f(x)$, and let $\tilde g(x)$ denote its noisy approximation. 
We denote $g_k:=g(x_k)$ and $\tilde g_k:=\tilde g(x_k)$. 
Given a search direction $\tpk$, a straightforward extension of the Armijo sufficient decrease condition \cite{mybook} reads
\[   
    \tilde f(x_k +\beta_k \tpk) - \tilde f(x_k) \leq c \beta_k \tilde g_k^T \tpk, \quad c \in (0, 1].
\]
This requires the comparison of noisy function evaluations
(case 1 in \S\ref{sec:intro}) and can lead to
poor performance or failure \cite{berahas2019derivative,oztoprak2021constrained,sun2022trust}.
To see this, suppose e.g. that $\tpk = - \tilde g_k$. Then, the right hand side is  always negative, but  due to the noisy nature of $\tilde f$, the left hand side can be positive even for a very small steplength, forcing the line search to decrease $\beta_k$ even more.

 One approach for circumventing these difficulties is to introduce a  margin $\epsilon_A(x_k)$ and to relax the Armijo condition as follows \cite{berahas2019derivative,berahas2019global,shi2022noise,sun2022trust},
\begin{equation}  \label{reline}
    \tilde f(x_k +\beta_k \tpk) \leq \tilde f(x_k) + c \beta_k \tilde g_k^T \tpk + 2\epsilon_A(x_k).
\end{equation}
A gradient projection method using a relaxed line search is given in Algorithm~\ref{alg:NoisyOPt}. It depends on a parameter $\alpha_0$ that determines the initial trial point in the line search. The importance of $\alpha_0$ will be discussed in subsequent sections}. 
In the algorithm,  $P_\Omega[\cdot]$ denotes the projection operator onto the feasible region $\Omega$. 
For the moment, we assume that $\epsilon_A(x_k)$ depends on  $\epsilon_f(x_k)$, and will elaborate on the exact nature of this relationship in the next subsection.

\begin{algorithm2e}[!ht] 
\SetAlgoLined
\textbf{Input:} Initial point {$x_0$}, constants $\rho\in(0,1)$, $c\in(0,1)$, 
and initial trial steplength  $\alpha_0>0$. \\
Set $k \leftarrow 0$. \\
\While{a termination condition is not met}{
Determine $\epsilon_A(x_k)$.\\
Compute a stochastic gradient $\tilde{g}_k$.\\
 $\tpk  \leftarrow P_{\Omega}[x_k- \alpha_0\tilde g_k]-x_k$.\\
 Set $\beta_k\leftarrow 1$. \\
 \While{$\tilde f({x}_k+ \beta_k \tpk)>\tilde f(x_k)+c\beta_k \tilde g_k^T\tpk+2\epsilon_A(x_k)$}{
 $\beta_k\leftarrow \rho \beta_k$.\\
 }
 $x_{k+1}\leftarrow x_k + \beta_k \tpk$. \\
 Set $k\leftarrow k+1$.
 }
\caption{({\tt GP-LS})  Line Search Gradient Projection Method }
\label{alg:NoisyOPt}
\end{algorithm2e}
In our experiments we use the parameters $\rho=1/2$ and $c=10^{-4}$.
We could have considered a more sophisticated gradient projection method with a projected backtracking line search \cite{bertsekas2015convex}, but the numerical and theoretical results would not be significantly different.

We now discuss the unspecified aspects of Algorithm~\ref{alg:NoisyOPt}, namely the computations of the relaxation $\epsilon_A(x_k)$ and the noisy gradients $\tgk$.

\subsection{\boldmath Choosing the Relaxation \texorpdfstring{$\epsilon_A(x)$}{epsA}} 
\label{sec:choosing_epsA}
One option is to choose $\epsilon_A(x)$ to be greater than $\epsilon_b$, where the latter is defined in \eqref{ebound} as a bound on the noise.  Then  \eqref{reline} is satisfied for all sufficiently small $\beta_k$, and 
one can establish deterministic convergence results to a neighborhood of the solution
\cite{berahas2019derivative,oztoprak2021constrained,sun2022trust}.
 However, in many applications, computing the bound $\epsilon_b$ is not feasible. Even when it is possible, choosing
$\epsilon_A(x) > \epsilon_b$ tends to be excessively cautious and 
 can degrade performance, as we will demonstrate in \S\ref{sec:numerical}.
 
 A more effective approach, in general, is to choose
$\epsilon_A(x)\leftarrow \lambda \epsilon_f(x)$, 
where $\epsilon_f(x)$ is the estimated noise level at $x$ and $\lambda$ is a positive constant. 
{We mentioned in Section~\ref{sec:noise} that estimating the bias is not required in designing the algorithm, for the following reasons.}

Suppose that the random variable $\Delta(x)$ is i.i.d. for all $x\in\Omega$, and that $\sigma_f(x)$ remains constant, so that computing $\epsilon_f(x)$ at a single $x$ suffices.
Then by utilizing concentration inequalities {(e.g. Chebyshev's or Vysochanskij–Petunin)} we can see that, $\mathbb{E}(\Delta(x))+\lambda \epsilon_f(x)$ serves as a high-probability estimate of $\epsilon_b$ for $\lambda$ large enough.
Given that the critical operations discussed in this paper solely involve comparisons or differences of function values, the {bias term $\mathbb{E}(\Delta(x))$} cancels out, justifying the rule $\epsilon_A(x)\leftarrow\lambda \epsilon_f(x)$.

This rule can also be motivated in the absence of the i.i.d assumption by introducing {the weaker set of assumptions: $\mathbb{V}(\Delta(x))\leq \sigma^2$ and $\mathbb{E}(\Delta(x))=0$ \cite{katya2018storm}.}
In that case it is reasonable to set $\epsilon_A(x)\leftarrow\lambda\sigma$.
Another line of research \cite{jin2021stepsearch,jin2023sample} that also motivates the rule $\epsilon_A(x)\leftarrow \lambda \epsilon_f(x)$ assumes the existence of probabilistic bounds of $\|\Delta(x)\|$, and allows for $\mathbb{E}(\Delta(x))\neq 0$. 

When the noise level does not vary significantly within the region of interest,  it is more
efficient to compute a constant estimate $\epsilon_f$ (as discussed  in \S\ref{sec:noise}) 
and fine tune the parameter $\lambda \in  [1, 2]$ to the application at hand. 
We can then drop the dependency on $x$ and write
\begin{equation}
 \label{eq:epsA_choice}
    \epsilon_A\leftarrow \lambda\epsilon_f.
\end{equation}
{In Section~\ref{sec:heuristics}, we will introduce the self-calibrated strategy which can adaptively tune $\epsilon_A$ and is efficient even with varying $\mathbb{E}(\Delta(x))$.}

In case the distribution of $\Delta(x)$ varies dramatically for different $x$, one may have to recompute $\epsilon_A$ during the course of the optimization or employ $\epsilon_b$ in lieu of a fixed value $\epsilon_A$; see Appendix~\ref{sec:limitations_of_noise} for details.

\subsection{Finite Difference Gradient Approximation}
\label{sec:FD}
The gradient of the objective function can be approximated  using (noisy) finite differences. 
This involves the critical operation 2 mentioned in \S\ref{sec:intro}.
To achieve stability, the function evaluations must be spread out appropriately to balance  truncation error and noise.

Let us consider the case where a universal noise level estimate  $\epsilon_f$ is available.
A value of $h$ that minimizes mean squared error for the forward difference estimator 
$$
    [\tilde g^{FD}(x)]_i:=\frac{\tilde f(x+he_i)-\tilde f(x)}{h},\quad i=1,\ldots,n ,
$$
is given by \cite{more2012estimating}
\begin{equation}
\label{eq:FDOptimalMSE}
    h\approx 8^{1/4}\sqrt{\frac{\epsilon_f}{L}},
\end{equation}
where $L$ is a bound on the second derivative of the objective function (or the Lipschitz constant of the gradient). (In this formula, $\epsilon_f$ should be replaced by $\epsilon_b$ when the latter is the only information available.)
Traditionally, the value of $L$ is estimated independently from $\epsilon_f$ \cite{GillMurrSaunWrig83,more2012estimating,shi2022numerical}. However,  Shi et al. \cite{shi2022adaptive} recently introduced a bisection procedure that calculates $h$ directly using only noisy evaluations $\tilde{f}$, avoiding a separate estimation of $L$.

In certain applications, such as the acoustic design problem described in \S\ref{sec:case}, analytic expressions for the gradient of a sample average approximation of the objective function are available; see Appendix~\ref{sec:horn_gradient}.
This will allow us to present a comparative efficiency analysis of noisy finite difference methods versus analytic gradients.

\section{Numerical Experiments} \label{sec:numerical}

We now describe  numerical experiments that test the efficiency of algorithms for solving noisy bound constrained optimization problem under various noise regimes. 
We  compare the line search gradient projection method {\tt GP-LS} defined in Algorithm~\ref{alg:NoisyOPt} with a variant using a fixed steplength, referred to as \texttt{GP-F} {(abbreviation of Gradient Projection Fixed-steplength Method)},  given by
\begin{equation}
    x_{k+1} \leftarrow P_{\Omega}[x_k- \alpha \tilde g_k] ,
\end{equation}
where $\tilde g_k$ is a gradient approximation, $P_\Omega[\cdot]$ is the projection operator onto the feasible region,
and $\alpha$ is a fixed steplength determined at the start of the algorithm.

Unless otherwise noted, the algorithms tested in this paper operate in the sample inconsistent case, meaning that every evaluation of the function uses a different batch of samples.
This applies both to finite difference approximations of gradients and to line searches.
(As a benchmark, we  report the results for the sample consistent case in Appendix~\ref{sec:sample_consistency}.)

\subsection{ Relaxed Line Search vs. Fixed Step Lengths}
\label{sec:LS_vs_SGD}
It is common practice to avoid line searches when minimizing noisy functions. 
We investigate whether this practice is still justified when employing the relaxed line search \eqref{reline}. To do so, we test our acoustic design problem under increasing noise levels.

In the first set of experiments, we compare the two gradient projection algorithms, {\tt GP-F} and {\tt GP-LS}, using gradients generated by finite differences. 
We chose a sample size $N=100$ in \eqref{precisely} for which the estimated noise level $\epsilon_f(b)$ varies between $10^{-3}$ and $10^{-2}$ (see Figure~\ref{fig:std100}).
Since $\epsilon_f(b)$ does not change dramatically, we use a single value $\epsilon_f$. We set $\epsilon_A =10^{-3}$ through \eqref{eq:epsA_choice}, after experimenting with the value of $\lambda$. 
Similarly, we use a fixed finite difference interval $h=10^{-2}$  in both methods, based on formula~\eqref{eq:FDOptimalMSE}  (experiments for other values of $h$ are discussed in the next subsection).
 
\begin{figure}[!htbp]
    \centering
    \includegraphics[width = 0.49\textwidth]{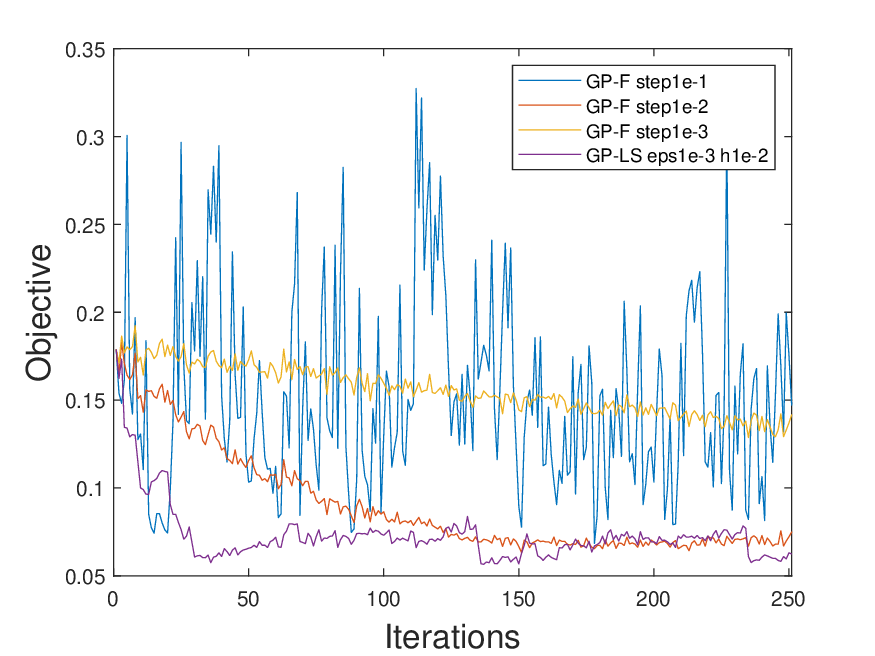}
    \includegraphics[width = 0.49\textwidth]{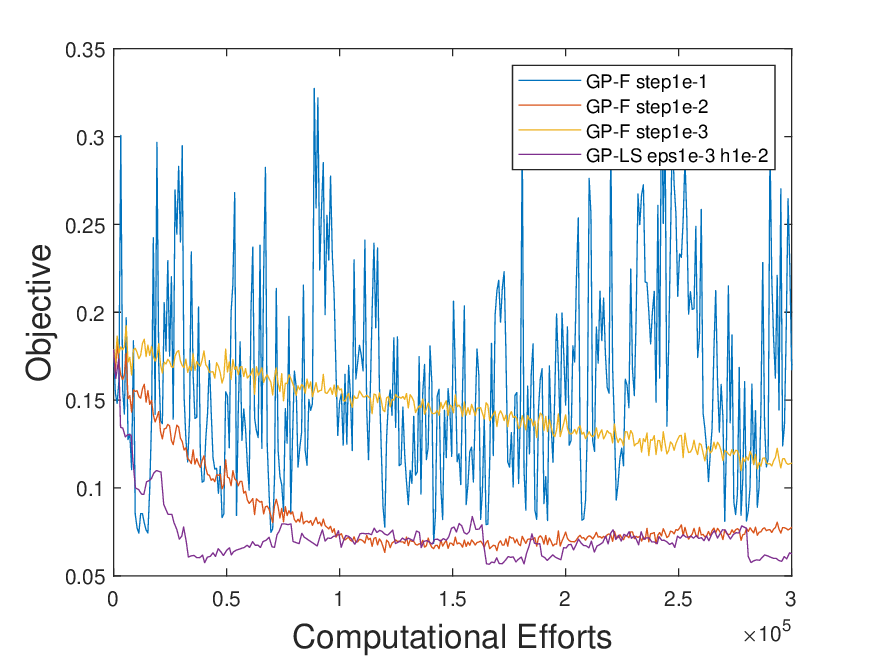}
    \caption{Comparison of the gradient projection method with ({\tt GP-LS}) and without ({\tt GP-F}) a line search;  the former using a relaxation  $\epsilon_A= 10^{-3}$ and the latter  using three values of $\alpha$. All methods use $N=100$ and a finite difference interval $h = 10^{-2}$. Left: Objective function value vs. iteration. Right: Objective function value vs. computational effort.}
    \label{fig:ls_vs_fixed}
\end{figure}

 The results are displayed in Figure~\ref{fig:ls_vs_fixed}.  Algorithm {\tt GP-F} was tested using three values of the fixed steplength, $\alpha = 10^{-1}, 10^{-2}, 10^{-3}$.  Algorithm {\tt GP-LS} 
 used an initial trial steplength  $\alpha_0=1$. 
 In the vertical axis we plot an approximation of the true objective function obtained by setting $N=100$ in \eqref{precisely}.
 In the left panel, the horizontal axis plots the iteration number;
 and in the right panel, it plots computational effort, defined as 
 \begin{equation}
 \label{comp_efforts_def}
     N\times \text{number of function calls}.
 \end{equation}

We observe from Figure~\ref{fig:ls_vs_fixed} that the performance of   {\tt GP-F} varies greatly with the choice of steplength $\alpha$. The value $\alpha=10^{-3}$ leads to a slow method, whereas the choice $\alpha=10^{-1}$ results in wild oscillations. The best performing method, using $\alpha= 10^{-2}$, was identified after extensive experimentation.
Observe that  {\tt GP-LS}  outperforms  the best option of {\tt GP-F}  in the initial third of the run.

\begin{figure}[!htbp]
    \centering
    \includegraphics[width = 0.49\textwidth]{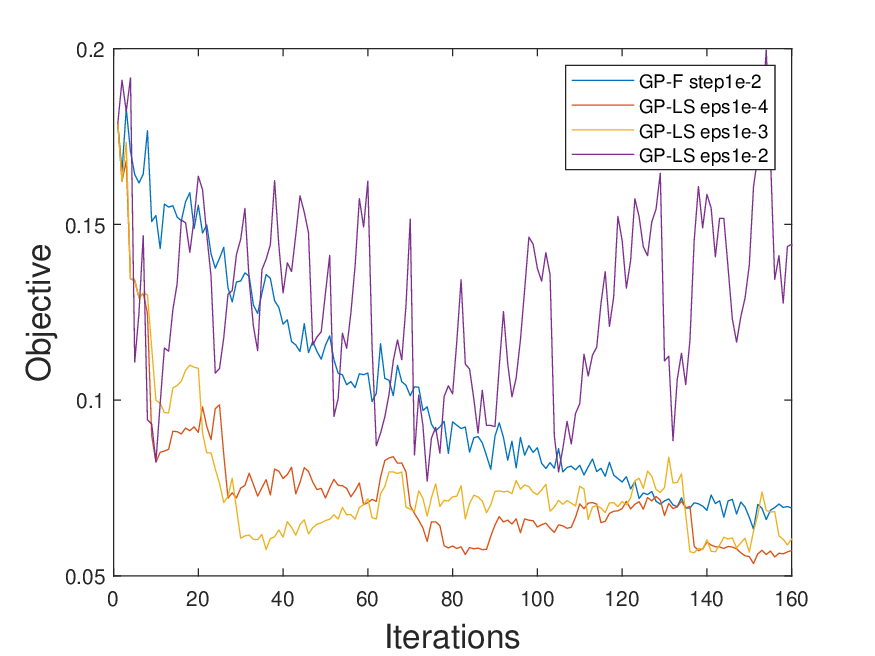}
    \includegraphics[width = 0.49\textwidth]{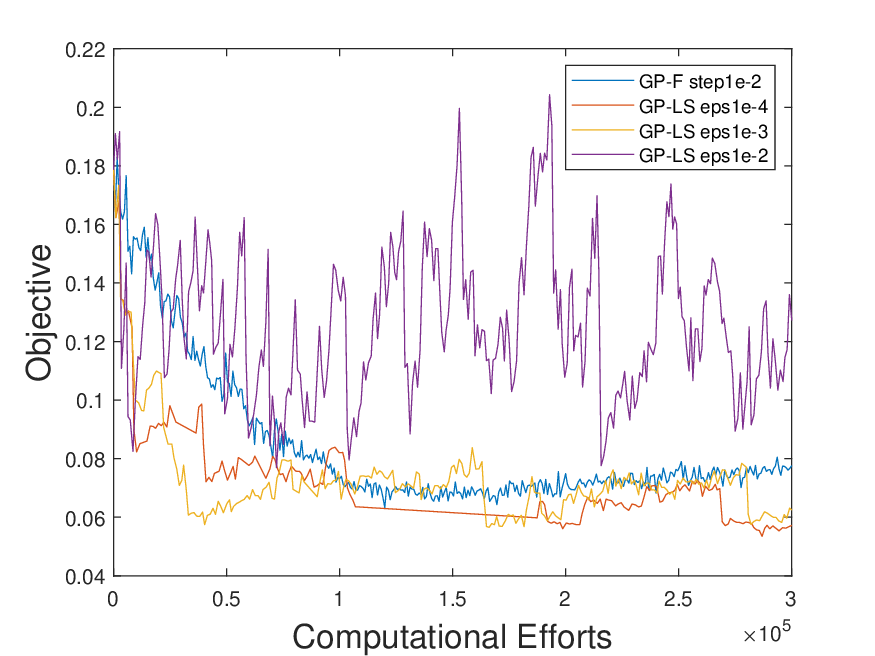}
    \caption{Performance of Algorithm {\tt GP-LS} with three values ($10^{-2}$, $10^{-3}$, $10^{-4}$) of the relaxation parameter $\epsilon_A$  in the line search.  We also plot the performance of Algorithm {\tt GP-F} with $\alpha = 10^{-2}$.
    Left: Objective function value vs. iteration. Right: Objective function  vs. computational effort \eqref{comp_efforts_def}.}
    \label{fig:ls_relax_comp}
\end{figure}

In the second set of experiments, we measure the effect of the relaxation parameter $\epsilon_A$ on algorithm~{\tt GP-LS}.   Figure~\ref{fig:ls_relax_comp} reports  results for choices $\epsilon_A= 10^{-2}, 10^{-3}, 10^{-4}$, which were derived as follows. 
For $N=100$, letting $\lambda=2$ and defining ${\epsilon}_f$ by \eqref{eq:epsf_avg}, we have that $\epsilon_A\approx 10^{-2}$ (such estimate is close to $\epsilon_b$).  To seek a lower bound {of $\epsilon_f(b)$},
we set $\lambda=1$, compute  $\epsilon_f$ by randomly sampling $b_1,\ldots,b_{100}$ in $\Omega$, 
and let $\epsilon_f=\min_{i=1,\ldots,100}\epsilon_f(b_i)$; this gives $\epsilon_A\approx10^{-3}$.
(We experiment with $\epsilon_A=10^{-4}$ in order to observe the effect of underestimating $\epsilon_A$.)

We observe from Figure~\ref{fig:ls_relax_comp} that {\tt GP-LS} performs well for $\epsilon_A= 10^{-3}$ and $ 10^{-4}$ but not so for $\epsilon_A=10^{-2}$. 
 By using this upper bound, the algorithm accepts overly noisy steps, resulting in oscillations. In contrast, if the relaxation $\epsilon_A$ is chosen too small (i.e., $10^{-4}$), it may cause the algorithm to repeatedly reject steps once it reaches the attainable accuracy in the function (observe the straight line in the right panel). 
However, this is not really harmful and a high number of rejections can be avoided by imposing a maximum number of backtracks; see e.g. the strategy in \S\ref{sec:heuristics}.
In summary, it is advisable to choose $\epsilon_A$ to be in the lower range of the estimated values of $\epsilon_f(b)$.

\subsection{ Finite Differences vs. Analytic Gradients}
\label{sec:numerical_FD_vs_analytic}
A common view in optimization is that finite difference gradient approximations  should be avoided in the noisy setting. We investigate this perspective in the context of the acoustic horn problem by comparing the use of finite differences and analytic expressions for the gradient {of a sample average approximation of the function. These analytic expressions are provided by the PDE solver as discussed in Appendix~\ref{sec:horn_gradient}.} 

In Figure~\ref{fig:fd_vs_ag}, we report the performance of the line search algorithm {\tt GP-LS} using finite differences or analytic gradients. 
We set $N=100$ and  $\epsilon_A=10^{-3}$ and {obtain the estimate} $h\approx 10^{-2}$ by using formula~\eqref{eq:FDOptimalMSE} with $\epsilon_f=10^{-3}$. Consequently, we report results with three values of the finite difference interval, namely $h= 10^{-1}, 10^{-2}, 10^{-3}$, to compare the outcomes of overestimating and underestimating interval choices.
In the figure on the right we plot the objective function vs. CPU time, which is an appropriate measure since the cost of an analytic gradient evaluation is difficult to quantify in terms of function evaluations.

The plots in Figure~\ref{fig:fd_vs_ag} indicate that, as anticipated,  the use of analytic gradients yields the best results. However,  the margin of improvement is not significant compared to {\tt GP-LS} with $h=10^{-2}$, a value of $h$ aligned with  formula~\eqref{eq:FDOptimalMSE}.

\begin{figure}[!htbp]
    \centering
    \includegraphics[width = 0.49\textwidth]{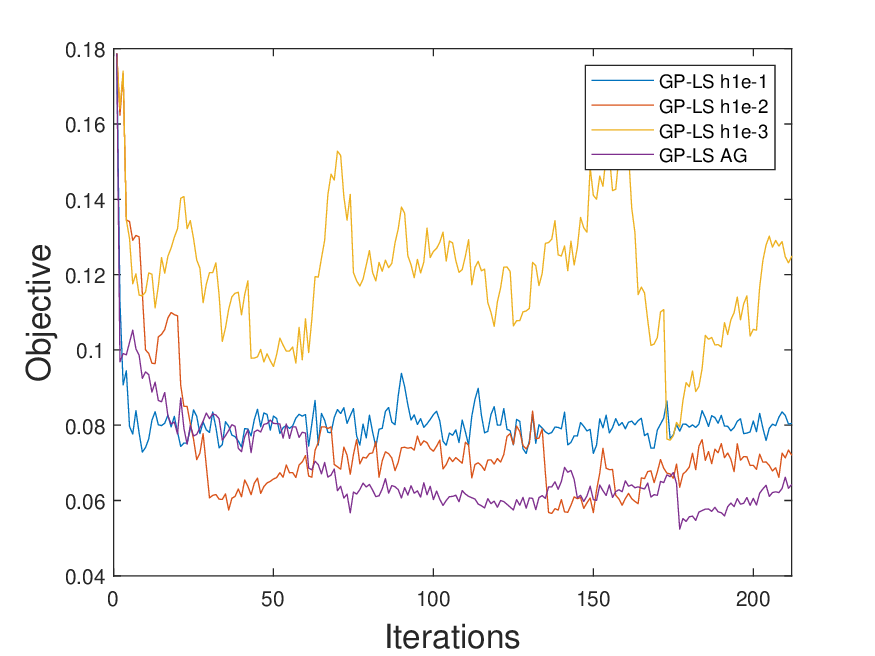}
    \includegraphics[width = 0.49\textwidth]{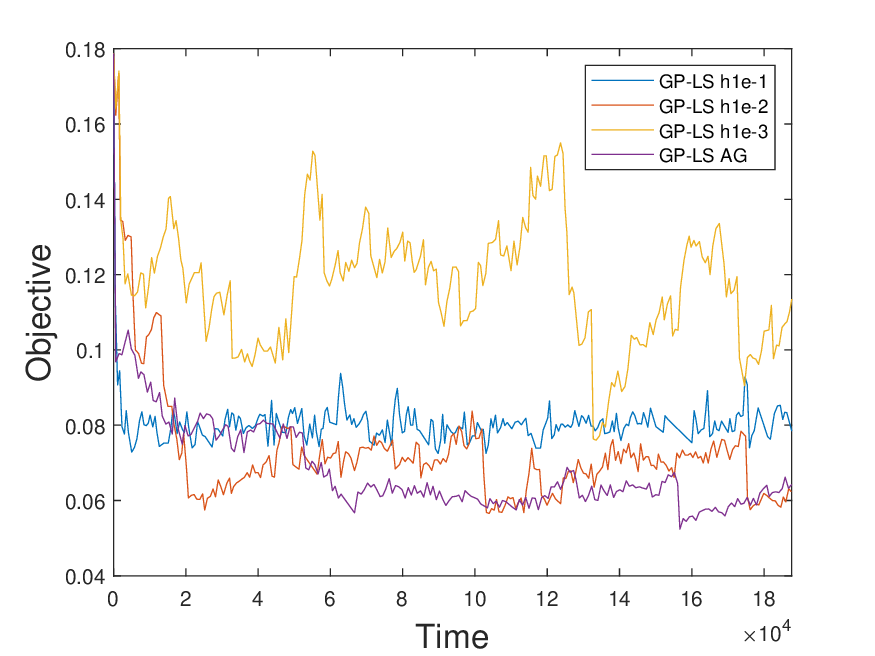}
    \caption{Comparison of analytic vs. finite difference  gradients in Algorithm~{\tt GP-LS}.
    We report results for three values of the finite difference parameter $h= 10^{-1}, 10^{-2}, 10^{-3}$. Left: Objective function value vs. iteration. Right: Objective function value vs.  CPU time.}
    \label{fig:fd_vs_ag}
\end{figure}

\subsection{\boldmath Increasing the Noise Level: \texorpdfstring{$N=50$}{N50}}  
\label{sec:nsamp50}
As the sample size decreases from $N=100$ to $N=50$, the problem becomes more noisy, potentially compromising the stability of the line search. Now, the convergence theory of stochastic gradient methods \cite{bottou2018optimization} states that the steplength should diminish in response to rising noise. This fact can be used to make the line search more robust by
decreasing the initial trial steplength $\alpha_0 $ in {\tt GP-LS}. 

\begin{figure}[!htbp]
    \centering
    \includegraphics[width = 0.49\textwidth]{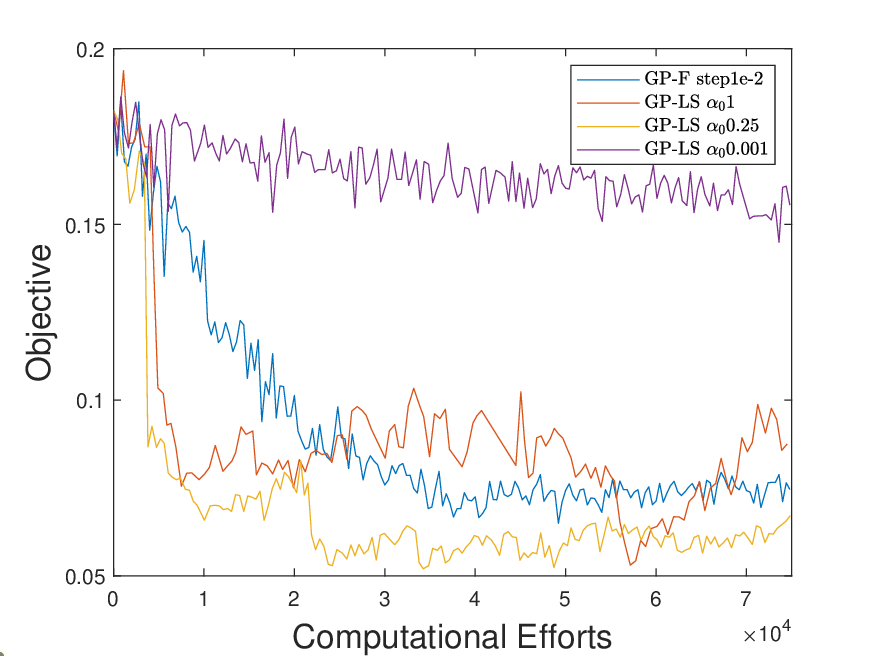}
    \includegraphics[width = 0.49\textwidth]{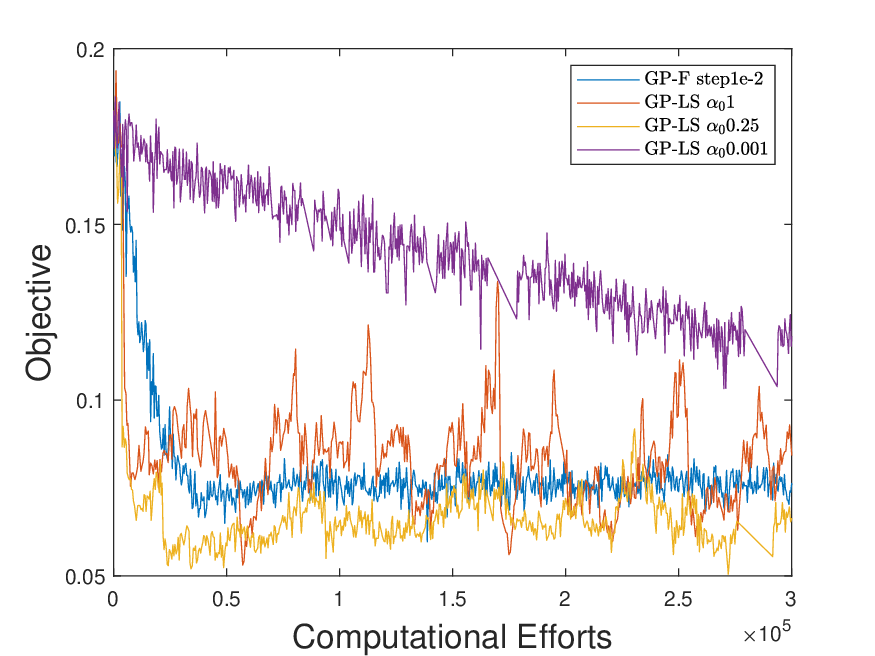}
    \caption{Comparison of three different values of the initial trial steplength, namely $\alpha_0 = 1, 0.25, 10^{-3}$ in {\tt GP-LS}. We also report {\tt GP-F} with $\alpha=10^{-2}$. Both algorithms were tested using $N=50$. Left: Objective function value vs. computational effort (up to $75,000$). Right: Objective function value vs. computational effort (up to $3\times 10^5$).}
    \label{fig:nsamp50}
\end{figure}

In Figure~\ref{fig:nsamp50}, we set $N=50$, $\epsilon_A=2\times 10^{-3}$, and plot the results for  {\tt GP-LS} with $\alpha_0 = 1,\, 0.25, 10^{-3}$. 
While $0.25$ is a reasonable choice, 1 and $10^{-3}$ are included to demonstrate the effects of excessively large or small choices of $\alpha_0$.
We also report the performance of 
{\tt GP-F} with $\alpha=10^{-2}$, a steplength obtained via tuning. 
The two figures report objective value vs. computational effort (defined in~\eqref{comp_efforts_def}). 
The left panel focuses on the early stage of the run 
while the right panel plots the overall long term behavior.

Figure~\ref{fig:nsamp50} shows the benefits of using values of $\alpha_0$ smaller than 1  in {\tt GP-LS}. The choice $\alpha_0=0.25$ outperforms all other options including the tuned {\tt GP-F}. The very small value $\alpha_0=10^{-3}$  
leads to poor performance both because it limits the lengths of the steps unduly and because comparisons in the line search become unreliable, sometimes yielding repeated rejections of trial steplengths. 

\subsection{\boldmath A Higher Noise Level: \texorpdfstring{$N=10$}{N10}}
\label{sec:nsamp10}
\begin{figure}[!htbp]
    \centering
    \includegraphics[width = 0.49\textwidth]{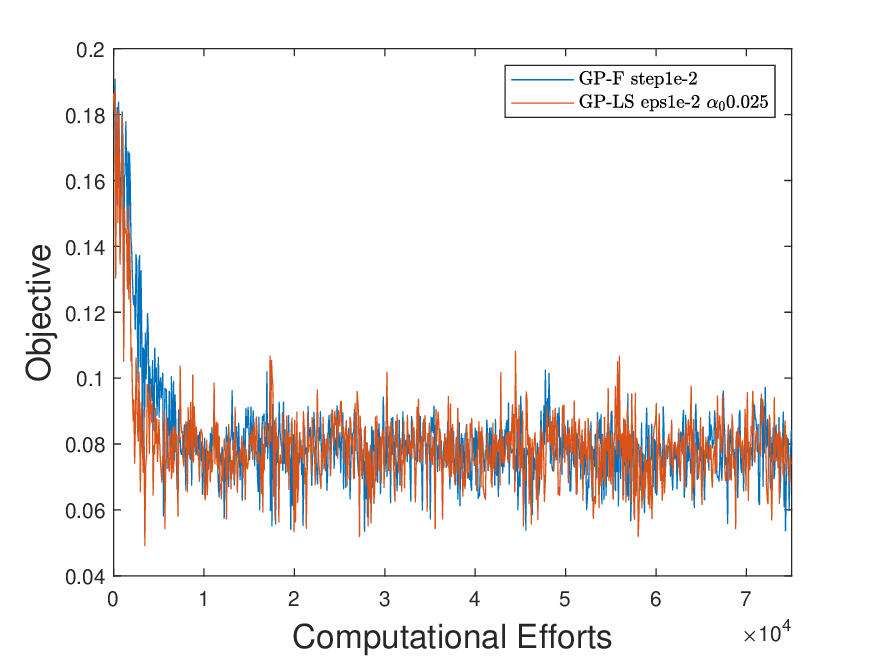}
    \includegraphics[width = 0.49\textwidth]{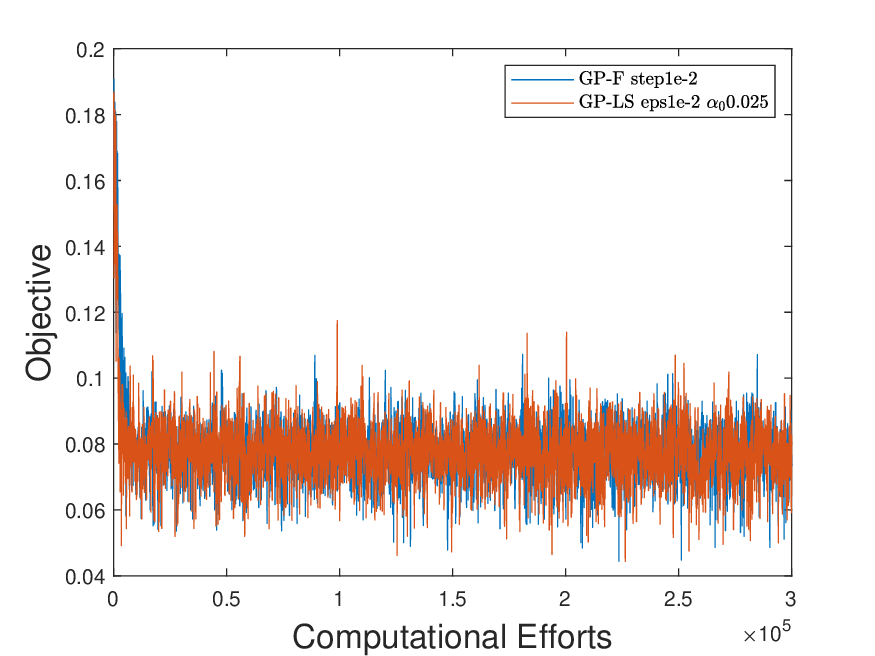}
    \caption{Comparison of {\tt GP-F} and {\tt GP-LS} with heuristics. Algorithm {\tt GP-LS} and {\tt GP-F} were tested using $N=10$; Left: Objective function value vs. computational effort (up to $75,000$). Right: Objective function value vs. computational effort (up to $3\times 10^5$).}
    \label{fig:nsamp10}
\end{figure}

When $N=10$, the noise level is so high that all algorithms exhibit strong oscillations in the objective.
In Figure~\ref{fig:nsamp10}, we report results of {\tt GP-F} with $\alpha=10^{-2}$, and {\tt GP-LS} with $\epsilon_A=10^{-2}$ and $\alpha_0=0.025$ (all parameters  chosen after  experimentation). The panel on the left focuses on the initial stages of the run, and the right panel on the overall run. Note that the best objective value achievable by the methods is around $8 \times 10^{-2}$, whereas for $N=50,100$ it was $6\times 10^{-2}$.
 {\tt GP-LS} no longer has an advantage over {\tt GP-F}, unlike the case for $N=50$ or $100$. 

To summarize our experiments so far, the relaxed line search strategy performs efficiently in the presence of noise by reducing the initial trial point $\alpha_0$ as the noise level increases. Yet, when dealing with highly noisy functions, employing a fixed step length strategy is  equally effective. Nevertheless, we now demonstrate that further enhancements to the line search strategy are possible.

\subsection{A Self-Calibrated Line Search Strategy}
\label{sec:heuristics}

We now show that the performance of the {\tt GP-LS} method can be improved significantly in the highly noisy regime by adaptively selecting the two key parameters in the {\tt GP-LS} method: $\epsilon_A$ and $\alpha_0$.
To do so, we first define a user-specified memory size $T$.
Every $T$ iterations, before computing the noisy gradients in \texttt{GP-LS}, instead of performing line 4 of Algorithm~\ref{alg:NoisyOPt}, we proceed as follows:
\begin{itemize}
    \item Compute the average number of line search backtracks in the most recent $T$ iterations, denoted as $avg$.
    \item If $avg\geq 3$, then update $\epsilon_A$ and $\alpha_0$ as
    \begin{equation}  \label{pat1}
    \epsilon_A\leftarrow \min\{1.5\epsilon_A,2\epsilon_f\},\quad \alpha_0\leftarrow 0.5\alpha_0,
    \end{equation}
    and if $avg\leq 0.1$,
    \begin{equation}  \label{pat2}
        \epsilon_A\leftarrow 0.5\epsilon_A, \quad \alpha_0\leftarrow \min\{1.5\alpha_0,10^{-1}\}.
    \end{equation}
\end{itemize}
The motivation for this strategy is as follows.

{Case 1:} If $avg$ is large, then either the relaxation is too small and the line search has stagnated (see Figure~\ref{fig:ls_relax_comp} for $\epsilon_A=10^{-4}$),
or the search direction is too noisy leading to many backtracking steps.
In this case,  the strategy increases $\epsilon_A$ and decreases $\alpha_0$ to further relax the line search and put more emphasis on safeguarding  errors. 
The upper bound of $\epsilon_A$ is set as $2\epsilon_f$ since we have seen in \S\ref{sec:choosing_epsA} that line search will ultimately be successful with high probability as $\epsilon_A$ is increased to $2\epsilon_f$.

{Case 2:} If $avg$ is small, then either $\epsilon_A$ is adequately large or the steps are productive.
In this case, we decrease $\epsilon_A$ since we try to keep this parameter as small as possible, and increase $\alpha_0$ to attempt to take more aggressive steps.

In addition to the rules \eqref{pat1} and \eqref{pat2}, we limit the number of possible backtracks by requiring that $\beta_k$ never be smaller than $\rho^{3T}$, where $\rho$ is the contraction parameter defined in Algorithm~\ref{alg:NoisyOPt}. Thus, the condition in the while loop in line 8 of Algorithm~\ref{alg:NoisyOPt} is changed to
\begin{equation}
    \tilde f({x}_k+ \beta_k \tpk)>\tilde f(x_k)+c\beta_k \tilde g_k^T\tpk+2\epsilon_A \text{ and }
    \beta_k \geq \rho^{3T}.
\end{equation}
Whenever $\beta_k$ is smaller than $\rho^{3T}$, the current step $\tpk$ will be discarded.


The constants in \eqref{pat1} and \eqref{pat2} can be tuned for the application at hand, but the method is not sensitive to the choices of these constants, with one caveat. 
It is important that, when changing $\epsilon_A$ and $\alpha_0$, we decrease them more rapidly than increase them (note $1.5\times 0.5<1$) because it is less harmful to perform more backtracks than accepting a poor step.
{We mention in passing that this method stands in contrast to a recently proposed method \cite{sun2023stochastic}, where an estimation of the gradient norm variance was used to re-scale the steps.}

\begin{figure}[!htbp]
    \centering
    \includegraphics[width = 0.49\textwidth]{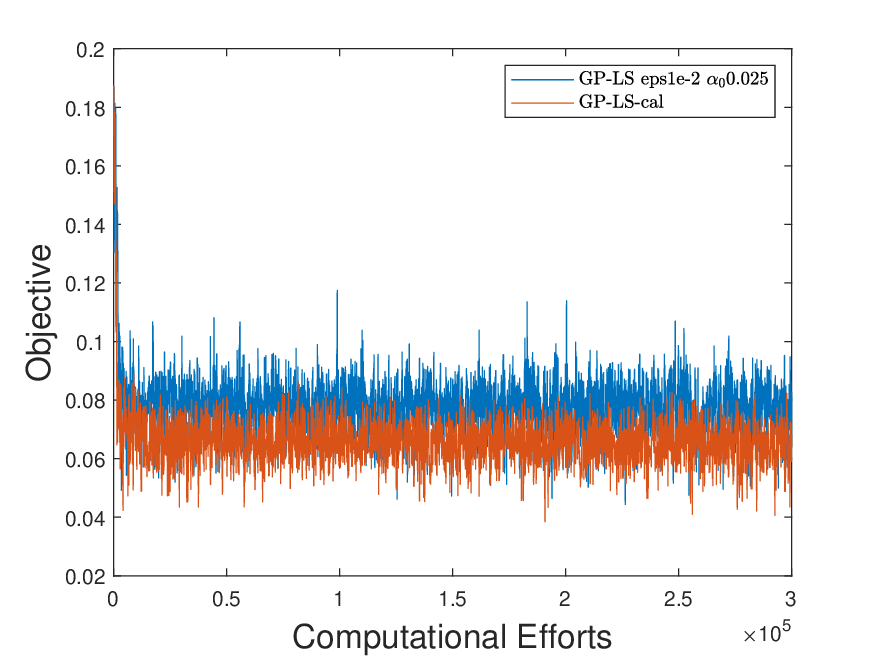}
    \includegraphics[width = 0.49\textwidth]{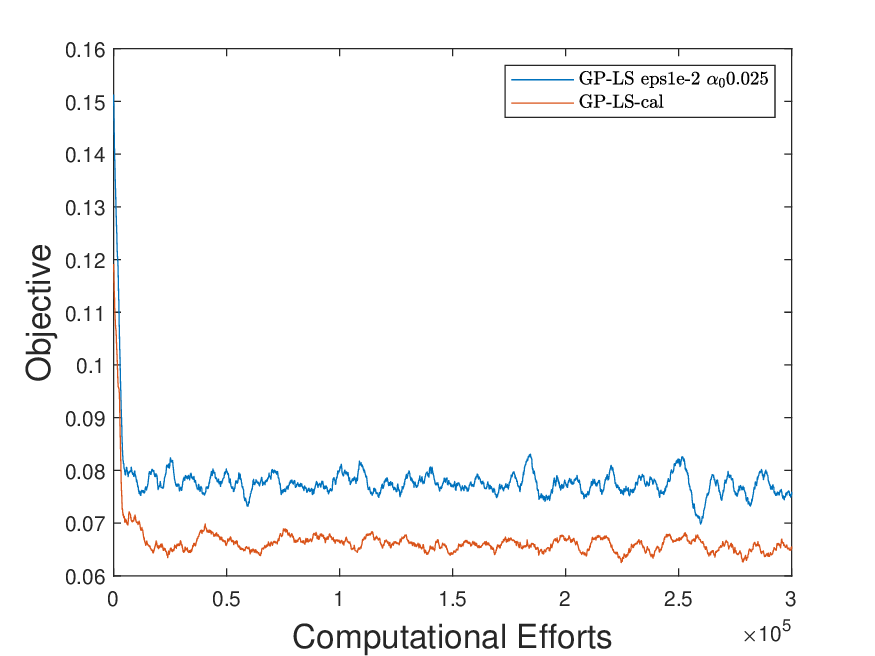}
    \caption{Comparison of {\tt GP-LS} and {\tt GP-LS-cal}. Algorithms were tested using $N=10$, $T=5$; Left: Objective function value vs. computational effort. Right:  Moving average of recent 50 objective function values vs. computational effort.}
    \label{fig:nsamp10_strategy}
\end{figure}

The results of applying $\texttt{GP-LS}$ with the self-calibrated strategy, denoted as \\ $\texttt{GP-LS-cal}$, are displayed in Figure~\ref{fig:nsamp10_strategy}.  There,  $T=5$ and the sample size is $N=10$.
The left panel compares  fine-tuned $\texttt{GP-LS}$ against $\texttt{GP-LS-cal}$. The right panel plots a smoothed version of the left figure, i.e., a moving average of objective.
We can observe that $\texttt{GP-LS-cal}$ clearly outperforms {\tt GP-LS}.
Moreover, the best average objective value of {\tt GP-LS-cal} improves to around $6\times 10^{-2}$, which is similar to the objective obtained for $N=50$ and $100$.

\section{Convergence Analysis}
\label{sec:convergence}
In this section, we establish convergence properties for algorithm {\tt GP-LS} in the presence of bounded noise, when applied to the problem
\begin{equation}
\label{opt_problem}
    \min_{x\in\Omega}\quad f(x),
\end{equation}
where $f$ is a nonlinear function and $\Omega$ is a closed convex set. 
We begin by stating two common assumptions. 

\begin{assumption}
\label{assump:closed}
  $\Omega$ is a nonempty, closed, and convex set, and for any $x\in\Omega$, $f(x)>-\infty$. 
\end{assumption}
\begin{assumption}
\label{assump:lipschitz_grad}
    $f$ is continuously differentiable in the feasible region $\Omega$ and for all $x,y\in\Omega$, there exist $L>0$ such that
    $$ \|g(x)-g(y)\|\leq L\|x-y\|  .$$
\end{assumption}
\noindent Next, we assume that the noise in the function and gradient is bounded.
\begin{assumption}
\label{assump:bounded_noise}
     For all $x\in\Omega$, there exists a constant $\epsilon_b>0$ such that
    $$\|\tilde f(x)-f(x)\|\leq \epsilon_b.$$
\end{assumption}
\begin{assumption}
\label{assump:bounded_grad_noise}
    For all $x\in\Omega$, there exists a constant $\epsilon_g>0$ such that
    $$\|\tilde g(x)-g(x)\|\leq \epsilon_g.$$
\end{assumption}

\medskip
Let us comment on the last two assumptions. 
In many engineering applications, including our acoustic design problem, the noise in the objective and gradient is inherently bounded due to the physical nature of the problem. In that case, when employing finite difference gradient estimations, we have that  Assumption~\ref{assump:lipschitz_grad} and~\ref{assump:bounded_noise} imply Assumption~\ref{assump:bounded_grad_noise}. This justifies the results in  \cite{berahas2019derivative,berahas2019global,berahas2021sequential,oztoprak2021constrained,sun2022trust} and the analysis given below, which assume bounded noise.

Nonetheless, a series of studies  
\cite{cao2023first,fang2022fully,jin2021stepsearch,jin2023sample}
assume only probabilistic bounds on the noise, and as mentioned in \S\ref{sec:choosing_epsA}, achieve high-probability convergence results. That analysis is more sophisticated but also more involved than the one presented here. Since we believe that the boundedness assumption holds in many  applications, our analysis is relevant to practice. It is also  novel in that no prior results exists for noisy gradient projection methods \emph{with a line search}, to our knowledge.

We begin the proof of global convergence by citing several established lemmas and introducing a stationarity measure specifically tailored to this problem.
Our ultimate objective is to demonstrate that the limit inferior of this measure is of order $O(\eps_b +\eps_g^2$).

\begin{lemma}[Prop. 1.1.9 Appendix B in \cite{bertsekas2015convex}]
\label{lem:uniqueness} 
    For any $x\in\mathbb{R}^n$, the projection of $x$ on $\Omega$ exists and is unique. 
    Furthermore, $z$ is the projection of $x$ on $\Omega$ if and only if $(x-z)^T(y-z)\leq 0$ for all $y\in \Omega$.
\end{lemma}

\begin{lemma}[Theorem 9.5-2 part (5) in \cite{kreyszig1991introductory}]
\label{lem:ProjContract}
    For any $x,y\in\mathbb{R}^n$, 
    $$ 
    \|P_\Omega[x]-P_\Omega[y]\|\leq \|x-y\|.$$
\end{lemma}

\noindent We now recall a standard stationary measure from optimization \cite{bertsekas2015convex}: 
\begin{equation}
         p(x):=P_\Omega[x-\alpha_0 g(x)]-x, \qquad
         \tilde{p}(x):=P_\Omega[x-\alpha_0\tilde{g}(x)]-x,  \\
\label{psdef}
\end{equation}
where $\alpha_0>0$ is any initial step-length set in Algorithm~\ref{alg:NoisyOPt}. {Note that by design of our gradient projection method, $p(x)$ is the search direction.}

\begin{lemma}
\label{lem:StatExist}
     $x^*\in \Omega$ is a first-order stationary point of problem \eqref{opt_problem} if and only if $p(x^*)=0$.
\end{lemma}
\noindent This lemma is a simple extension of a classical result (see Prop. 6.1.1 (b) in \cite{bertsekas2015convex});  we include its proof in Appendix~\ref{sec:supp_proof} for completeness.

\begin{remark}
\label{rem:stat}
Lemma~\ref{lem:StatExist} implies that beyond serving as the search direction for algorithm \texttt{GP-LS} at iteration $x_k$,  $p(x_k)$ also functions as a measure of stationary for problem~\eqref{opt_problem}. 
There is, however, another optimality measure that is more convenient in deriving our main convergence result. This measure is given by $-p(x_k)^Tg(x_k)$, as discussed next.
\end{remark}





\begin{lemma}
\label{lem:stat}
    $ x^*\in \Omega$ is a first-order stationary point of problem \eqref{opt_problem} if and only if 
    \begin{equation}  \label{newstat}
    p(x^*)^T g(x^*)=0.
     \end{equation}
\end{lemma}

\begin{proof}
    By Lemma~\ref{lem:StatExist}, it suffices to show that \eqref{newstat} is equivalent to $p(x^*)=~0$. Clearly $p(x^*)=0\ \Rightarrow\ p(x^*)^T g(x^*)=0$.     

To establish the converse, assume that $ p(x^*)^T g(x^*)=0$, and define $\theta$ as the angle between $p(x^*)$ and $ g(x^*)$, so that 
$$\|p(x^*)\|\| g(x^*)\|\cos\theta = 0. $$ 
If $\|p(x^*)\|=0$ or $\| g(x^*)\|=0$ (which by \eqref{psdef} implies $\|p(x^*)\|=0$), then $p(x^*) = 0$, yielding the desired result. 

Let us therefore consider the case when $\|p(x^*)\|\neq 0$ and $\| g(x^*)\|\neq0$, and $\cos\theta = 0$. We show by contradiction that this case is not possible. 
    Note from \eqref{psdef}
    \begin{equation*}
    \begin{split}
        \|P_\Omega[x^*-\alpha_0 g(x^*)]-(x^*-\alpha_0 g(x^*))\|^2 & =\|P_\Omega[x^*-\alpha_0 g(x^*)]-x^*\|^2+\|\alpha_0 g(x^*)\|^2 \\
        & > \|P_\Omega[x^*-\alpha_0 g(x^*)]-x^*\|^2.
    \end{split}
    \end{equation*} 
    This contradicts the fact that $P_\Omega[x^*-\alpha_0 g(x^*)]$ as the unique vector closest to $x^*-\alpha_0 g(x^*)$ in $\Omega$.
\end{proof}
{Using the standard abbreviations $p_k:=p(x_k)$, $\tpk:=\tilde p(x_k)$,
Lemma~\ref{lem:stat} establishes $p_k^Tg_k$ as a stationary measure of problem~\eqref{opt_problem}---and $\tpk^T\tgk$ is its noisy counterpart, which is the quantity accessed by the algorithm.
In light of Lemma~\ref{lem:uniqueness}, it is easy to see that $-p_k^Tg_k\geq 0$ and $-\tpk^T\tgk\geq 0$.}

Let us now define
\begin{equation}
    \delta_g(x_k):=(-\tilde g_k)-(-g_k),\qquad \delta_p(x_k):=\tilde{p}_k-p_k.
\end{equation}

We now establish a technical lemma relating
$-\tpk^T\tgk$ and the stationary measure $-p_k^Tg_k$, in terms of a scaling {factor} dependent on the magnitude of the noise $\|\delta_g(x)\|$.

\begin{lemma}
\label{lem:ineq_noisy_stat_measure}
    Under the assumptions previously stated, for any iterate $x_k$ generated by {\tt GP-LS} (Algorithm~\ref{alg:NoisyOPt}),
    $$-\tpk^T\tgk \geq \left( \frac{1}{2}-\frac{\alpha_0}{2}-\left(\frac{3\alpha_0}{2}+\frac{1}{2}\right )\gamma_k^2 \right)(-p_k^Tg_k)$$
    where 
    \begin{equation}
        \label{eq:gamma_k_def}
        \gamma_k:=\frac{\|\dgk\|}{\sqrt{-p_k^Tg_k}}.
    \end{equation}
\end{lemma}

\begin{proof}

We lead the proof by noting the differences between $-p_k^Tg_k$ and $-\tpg$:
     \be 
      -\tpg - (-p_k^Tg_k) = -g_k^T\dpk+p_k^T\dgk+\dpk^T\dgk.
     \label{eq:terms}
     \ee
     We establish bounds on each terms on the right hand side of this equation. 

     We first show that the last term $\dpk^T\dgk$ is non-negative.
    Apply Lemma~\ref{lem:uniqueness} with $x=x_k-\alpha_0g_k$, $z =  P_\Omega[x_k-\alpha_0 g_k] =  x_k+p_k $, and $y=x_k+\tpk$, we have $(-\alpha_0g_k-p_k)^T(\tpk-p_k)\leq 0$, which implies
    \begin{equation}
    \label{ineq:gdpkLB_no_noise}
        -p_k^T\dpk \leq \alpha_0g_k^T\dpk.
    \end{equation}
    Apply again Lemma~\ref{lem:uniqueness} with $x=x_k-\alpha_0\tilde{g}_k$,
    $ z = P_\Omega[x_k-\alpha_0\tilde{g}_k] = x_k+\tilde{p}_k$ 
    and $y=x_k+p_k \in\Omega$, we have $(-\alpha_0\tgk-\tpk)^T(p_k-\tpk)\leq 0$, which implies
    \begin{equation}
    \label{ineq:gdpkLB}
        \tpk^T\dpk\leq -\alpha_0\tgk^T\dpk.
    \end{equation}
    Summing up~\eqref{ineq:gdpkLB_no_noise} and~\eqref{ineq:gdpkLB}, we obtain
    \begin{align}
         & (\tpk-p_k)^T\dpk\leq \alpha_0(g_k-\tgk)^T\dpk \nonumber \\
        \Longrightarrow \quad & \dpk^T\dgk \geq \frac1{ \alpha_0}\|\dpk\|^2 \geq 0.
     \label{ineq:dpkdgk_nonnegative}
    \end{align}

    We next analyze the cross term $g_k^T\dpk$. For this, we first derive a few auxilary inequalities. First note by Lemma~\ref{lem:ProjContract}, 
    \begin{equation}
    \label{ineq:PnoiseBound}
        \|\delta_p(x_k)\|=\|p_k-\tilde{p}_k\|\leq \alpha_0\|g_k-\tilde{g}_k\|=\alpha_0\|\delta_g(x_k)\|.
    \end{equation}
    Moreover, apply Lemma~\ref{lem:uniqueness} with $x=x_k-\alpha_0g_k$, $z =  P_\Omega[x_k-\alpha_0 g_k] =  x_k+p_k $, and $y=x_k\in\Omega$, we obtain
    \begin{equation}
    \label{ineq:ProjLemNoiseFree}
        \|p_k\|^2\leq -\alpha_0p_k^Tg_k.
    \end{equation}
    To bound $g_k^T\dpk$, we have from~\eqref{ineq:gdpkLB} 
    \begin{equation}
        \tpk^T\dpk\leq -\alpha_0\tgk^T\dpk=-\alpha_0g_k^T\dpk+\alpha_0\dpk^T\dgk.
    \end{equation}
    Re-organize and obtain
    \begin{align}
        -\alpha_0g_k^T\dpk & \geq \tpk^T\dpk-\alpha_0\dpk^T\dgk \nonumber\\
        &  =   p_k^T\dpk+\dpk^T\dpk-\alpha_0\dpk^T\dgk \nonumber\\
        & \geq -\alpha_0\|\dgk\|\|p_k\|-\alpha_0^2\|\dgk\|^2 \nonumber\\
        & \geq -\alpha_0\|\dgk\|\sqrt{-\alpha_0p_k^Tg_k}-\alpha_0^2\|\dgk\|^2 \nonumber\\
        & \geq \frac{\alpha_0}{2}p_k^Tg_k-\frac{3\alpha_0^2}{2}\|\dgk\|^2 \nonumber\\
        &  = \left(\frac{\alpha_0}{2}+\frac{3\alpha_0^2}{2} \gamma_k^2\right) p_k^Tg_k.
    \label{ineq:-gkdpk_LB_gamma}
    \end{align}
    Here, the second inequality follows from Cauchy-Schwartz inequality, $\|\delta_p(x_k)\|^2\geq 0$, and \eqref{ineq:PnoiseBound}; 
    the third inequality follows from~\eqref{ineq:ProjLemNoiseFree}; the fourth is from arithmetic-geometric mean, i.e., $\alpha_0\|\dgk\| \sqrt{-\alpha_0p_k^Tg_k}\leq \frac{1}{2}\left( 
    \alpha_0^2\|\dgk\|^2-\alpha_0p_k^Tg_k\right)$; and the last line is by $\gamma_k$ defined in~\eqref{eq:gamma_k_def}.

    
    Finally, using~\eqref{ineq:-gkdpk_LB_gamma},~\eqref{ineq:dpkdgk_nonnegative} \& Cauchy-Schwartz,~\eqref{ineq:ProjLemNoiseFree}, and arithmetic-geometric mean for the following inequalities respectively, we obtain the desired result:
    \begin{align}
        -\tpg & = -p_k^Tg_k+(-g_k^T\dpk)+p_k^T\dgk+\dpk^T\dgk \nonumber \\
        & \geq -p_k^Tg_k+\left( \frac{1}{2}+\frac{3\alpha_0}{2}\gamma_k^2\right)p_k^Tg_k+p_k^T\dgk+\dpk^T\dgk \nonumber\\
        & \geq \left( \frac{1}{2}-\frac{3\alpha_0}{2}\gamma_k^2\right)(-p_k^Tg_k)-\|p_k\|\|\dgk\| \nonumber\\
        & \geq \left(\frac{1}{2}-\frac{3\alpha_0}{2}\gamma_k^2 \right)(-p_k^Tg_k)-\|\dgk\|\sqrt{-\alpha_0p_k^Tg_k} \nonumber\\
        & \geq \left(\frac{1}{2}-\frac{3\alpha_0}{2}\gamma_k^2 \right)(-p_k^Tg_k)-\frac{1}{2}\|\dgk\|^2+\frac{\alpha_0}{2}p_k^Tg_k \nonumber\\
        & =\left( \frac{1}{2}-\frac{\alpha_0}{2}-\left (\frac{3\alpha_0}{2}+\frac{1}{2} \right)\gamma_k^2 \right)(-p_k^Tg_k).
        \label{ineq:LBUB}
    \end{align}
\end{proof}

We can now state the main convergence theorem for the gradient projection algorithm with a relaxed line search. 
We recall that $-p_k^Tg_k$ serves both an algorithmic role in the Armijo decrease condition and a theoretical role as a stationary measure of the problem, as mentioned in Remark~\ref{rem:stat}.

\begin{theorem}
    Under Assumptions~\ref{assump:closed}-\ref{assump:bounded_grad_noise}, if $\alpha_0+2c<1$ and $\epsilon_A>\epsilon_b$, 
    the iterates $\{x_k\}$ generated by {\tt GP-LS} (Algorithm~\ref{alg:NoisyOPt}) satisfy
    $$\liminf_{k\to\infty} \left |p_k^Tg_k\right | \leq \Bar{\epsilon}$$
    where
    \begin{equation}
    \label{eq:def_bar_eps}
        \Bar{\epsilon}:=\frac{\epsilon_g^2}{\gamma^2}+\frac{2\alpha_0L}{c\rho\left( \frac{1}{2}-\frac{\alpha_0}{2}-\left(\frac{3\alpha_0}{2}+\frac{1}{2}\right)\gamma^2 \right)}\left(\epsilon_A+\epsilon_b\right),
        \end{equation}
        and
        \begin{equation}
        \label{eq:def_gamma_2}
        \gamma^2:=\frac{(1-2c-\alpha_0)(1-\alpha_0)}{(1-2c-\alpha_0)\left(3\alpha_0+1\right)+2}.
    \end{equation}
\end{theorem}

\begin{proof}

    The proof is constructed by characterizing the descent on the objective function using the noisy stationary measure $\tilde p_k^T \tilde g_k$, and dividing the proof into two cases according to the relative size of the noise.
    
    First, 
    by applying Lemma~\ref{lem:uniqueness} with $x=x_k-\alpha_0\tilde g_k$, $z =  P_\Omega[x_k-\alpha_0 \tilde g_k] =  x_k+\tilde p_k $, and $y=x_k\in\Omega$, we have
    \begin{equation}
    \label{ineq:tpkBound}
        \|\tilde{p}_k\|^2\leq -\alpha_0\tpg.
    \end{equation}
    Next, by a Taylor expansion and Assumption~\ref{assump:lipschitz_grad}, we have for any $\beta>0$,
    \begin{equation}
    \begin{split}
        f(x_k+\beta \tpk) & \leq f(x_k)+\beta\tpk^Tg_k+\frac{L}{2}\beta^2\|\tpk\|^2 \\
        & \leq f(x_k)+\beta \tpk^T(\tgk+\dgk)+\frac{L}{2}\beta^2(-\alpha_0\tpg) \\
        & \leq f(x_k)+(-\beta+\frac{\alpha_0L}{2}\beta^2)(-\tpg)+\beta\|\tpk\|\|\dgk\| \\
        & \leq f(x_k)+(-\beta+\frac{\alpha_0L}{2}\beta^2)(-\tpg)+\beta\|\dgk\|\sqrt{-\alpha_0\tpg} \\
        & \leq f(x_k)+\left(\left(\frac{\alpha_0}{2}-1\right)\beta+\frac{\alpha_0L}{2}\beta^2\right)(-\tpg)+\frac{\beta}{2}\|\dgk\|^2,
    \end{split}
    \end{equation}
    where the second and fourth inequalities are from~\eqref{ineq:tpkBound}, the third is from Cauchy-Schwartz, and the last is from the arithmetic-geometric mean.
    Together with Assumption~\ref{assump:bounded_noise}, we have
    \begin{equation}
    \label{ineq:RelaxedDescent}
        \tilde{f}(x_k+\beta \tpk)\leq \tilde{f}(x_k)+ \left[\left(\left(\frac{\alpha_0}{2}-1\right)\beta+\frac{\alpha_0L}{2}\beta^2\right)(-\tpg)+\frac{\beta}{2}\|\dgk\|^2 \right]+2\epsilon_b.
    \end{equation}
    
    We now note that the line search in \texttt{GP-LS} always terminates within finitely many back-tracking steps.
    This follows from the fact that we pick
     $\epsilon_A>\epsilon_b$ and  that the term inside square brackets in \eqref{ineq:RelaxedDescent} converges to zero as $\beta\to0$. Hence,
    the relaxed Armijo condition~\eqref{reline} will be satisfied for some sufficiently small $\beta_k>0$.

    We now divide the set of iterates into two cases depending on whether the noise dominates the optimality measure, in the sense that the ration  $\gamma_k$ is larger than the threshold $\gamma$, where these quantities are defined in \eqref{eq:gamma_k_def} and \eqref{eq:def_gamma_2}, respectively.

    Note by the assumption $\alpha_0+2c<1$ and simple algebra 
    \begin{equation}
    \label{ineq:gammaProp}
        0<\gamma^2<\frac{1-\alpha_0}{3\alpha_0+1}.
    \end{equation} 
    
    

    \textbf{Case 1}: Noise is relatively small: $\gamma_k^2\leq \gamma^2$. By~\eqref{ineq:RelaxedDescent}, \eqref{eq:gamma_k_def}, and Lemma~\ref{lem:ineq_noisy_stat_measure}, we have
    \begin{equation}
    \label{eq:case1_desc}
    \begin{split}
        \tilde{f}(x_k+\beta \tpk) & \leq \tilde{f}(x_k)+\left(\left(\frac{\alpha_0}{2}-1\right)\beta+\frac{\alpha_0L}{2}\beta^2\right)(-\tpg)+\frac{\beta}{2}\gamma_k^2(-p_k^Tg_k)+2\epsilon_b \\
        \leq  \tilde{f}(x_k) & +\left( \left(\frac{\gamma_k^2}{1-\alpha_0-(3\alpha_0+1)\gamma_k^2}+\frac{\alpha_0}{2}-1\right)\beta+\frac{\alpha_0L}{2}\beta^2\right)(-\tpg)+2\epsilon_b.
    \end{split}
    \end{equation}
    With this result, the Armijo condition \eqref{reline} holds when
    \be \left(\frac{\gamma_k^2}{1-\alpha_0-(3\alpha_0+1)\gamma_k^2}+\left(\frac{\alpha_0}{2}-1\right)\right)\beta+\frac{\alpha_0L}{2}\beta^2\leq -c\beta,\ee which is equivalent to
    \begin{equation}
    \label{ineq:betaUB}
        \beta\leq\frac{2}{\alpha_0L}\left( -c-\frac{\gamma_k^2}{1-\alpha_0-(3\alpha_0+1)\gamma_k^2}-\frac{\alpha_0}{2}+1 \right).
    \end{equation}
    Since $\gamma_k^2\leq \gamma^2$, by~\eqref{ineq:gammaProp},
    \begin{equation}
        1-\alpha_0-(3\alpha_0+1)\gamma_k^2 \geq 1-\alpha_0-(3\alpha_0+1)\gamma^2>0.
    \end{equation}
    With this, we note that $-\frac{\gamma_k^2}{1-\alpha_0-(3\alpha_0+1)\gamma_k^2}$ is decreasing in $\gamma_k^2$, for $\gamma_k^2\in(0,\gamma^2]$.
    Therefore its lower bound is achieved when $\gamma_k = \gamma$, i.e.
    \begin{equation}
        -c-\frac{\gamma_k^2}{1-\alpha_0-(3\alpha_0+1)\gamma_k^2}-\frac{\alpha_0}{2}+1\geq -c-\frac{\gamma^2}{1-\alpha_0-(3\alpha_0+1)\gamma^2}-\frac{\alpha_0}{2}+1=\frac{1}{2}
    \end{equation}
    where the equality follows from the definition of $\gamma^2$ in~\eqref{eq:def_gamma_2} and algebra. 
    
    This, together with~\eqref{ineq:betaUB}, implies that the relaxed Armijo condition holds for any $\beta\leq\frac{1}{\alpha_0L}$. Thus, for any $k\in\mathbb{N}$ in Case 1,
    \begin{equation}
    \label{ineq:beta_k_LB}
        \beta_k\geq\frac{\rho}{\alpha_0 L}.
    \end{equation}
    Therefore, we have from~\eqref{reline}, Assumption~\ref{assump:bounded_noise}, and~\eqref{ineq:beta_k_LB} that
    \begin{equation}
    \begin{split}
        f(x_k+\beta_k\tpk) & \leq f(x_k)+c\beta_k \tpg +2\epsilon_A+2\epsilon_b \\
        &\leq f(x_k)-\frac{c\rho}{\alpha_0L}\left( \frac{1}{2}-\frac{\alpha_0}{2}-\left(\frac{3\alpha_0}{2}+\frac{1}{2}\right)\gamma^2 \right)(-p_k^Tg_k)+2\epsilon_A+2\epsilon_b,
    \end{split} 
    \end{equation}
    which measures the reduction of the objective between any two consecutive iterations for Case 1 (notice that according to Algorithm~\ref{alg:NoisyOPt}, $x_{k+1}=x_k+\beta_k\tpk$).

    \textbf{Case 2}: Noise is relatively large: $\gamma_k^2>\gamma^2$. By definition of $\gamma_k$ in \eqref{eq:def_gamma_2}
    \begin{equation}
    \label{ineq:case2def}
        \|\dgk\|^2>\gamma^2(-p_k^Tg_k).
    \end{equation}
    As explained in the paragraph after \eqref{ineq:RelaxedDescent}, 
    there always exists $\beta_k>0$ such that the relaxed Armijo condition~\eqref{reline} holds. This fact, together with Assumption~\ref{assump:bounded_noise}, Assumption~\ref{assump:bounded_grad_noise}, and~\eqref{ineq:case2def},
    \begin{equation}  \label{light}
    \begin{split}
        f(x_k+\beta_k\tpk) \leq& f(x_k)-c\beta_k(-\tpg)+2\epsilon_A+2\epsilon_b \\
         \leq& f(x_k)+2\epsilon_A+2\epsilon_b \\
         =& f(x_k)-\frac{c\rho}{\alpha_0L}\left( \frac{1}{2}-\frac{\alpha_0}{2}-\left(\frac{3\alpha_0}{2}+\frac{1}{2}\right)\gamma^2 \right)(-p_k^Tg_k) \\ 
        & +\frac{c\rho}{\alpha_0L}\left( \frac{1}{2}-\frac{\alpha_0}{2}-\left(\frac{3\alpha_0}{2}+\frac{1}{2}\right)\gamma^2 \right)(-p_k^Tg_k)+2\epsilon_A+2\epsilon_b \\
        \leq& f(x_k)-\frac{c\rho}{\alpha_0L}\left( \frac{1}{2}-\frac{\alpha_0}{2}-\left(\frac{3\alpha_0}{2}+\frac{1}{2}\right)\gamma^2 \right)(-p_k^Tg_k) \\ 
        & +\frac{c\rho}{\alpha_0L}\left( \frac{1}{2}-\frac{\alpha_0}{2}-\left(\frac{3\alpha_0}{2}+\frac{1}{2}\right)\gamma^2 \right)\frac{\|\dgk\|^2}{\gamma^2}+2\epsilon_A+2\epsilon_b \\
         \leq &f(x_k)-\frac{c\rho}{\alpha_0L}\left( \frac{1}{2}-\frac{\alpha_0}{2}-\left(\frac{3\alpha_0}{2}+\frac{1}{2}\right)\gamma^2 \right)(-p_k^Tg_k) \\ 
        & +\frac{c\rho}{\alpha_0L}\left( \frac{1}{2}-\frac{\alpha_0}{2}-\left(\frac{3\alpha_0}{2}+\frac{1}{2}\right)\gamma^2 \right)\frac{\epsilon_g^2}{\gamma^2}+2\epsilon_A+2\epsilon_b \\
         = & f(x_k)-\frac{c\rho}{\alpha_0L}\left( \frac{1}{2}-\frac{\alpha_0}{2}-\left(\frac{3\alpha_0}{2}+\frac{1}{2}\right)\gamma^2 \right)(-p_k^Tg_k)+\eta,
    \end{split}
    \end{equation}
    where
    \begin{equation}
        \eta:=\frac{c\rho}{\alpha_0L}\left( \frac{1}{2}-\frac{\alpha_0}{2}-\left(\frac{3\alpha_0}{2}+\frac{1}{2}\right)\gamma^2 \right)\frac{\epsilon_g^2}{\gamma^2}+2\epsilon_A+2\epsilon_b.
    \end{equation}
    Condition \eqref{light} measures the reduction of the objective between any two consecutive iterations for Case 2.

    Now combine both Case 1 and 2, and since $\eta>2\epsilon_A+2\epsilon_b$, it follows that for all $k\in\mathbb{N}$,
    \begin{equation}
        f(x_{k+1})\leq f(x_k)-\frac{c\rho}{\alpha_0L}\left( \frac{1}{2}-\frac{\alpha_0}{2}-\left(\frac{3\alpha_0}{2}+\frac{1}{2}\right)\gamma^2 \right)(-p_k^Tg_k)+\eta.
    \end{equation}

    Finally, to prove that $\liminf_{k\to\infty}|p_k^Tg_k|\leq \bar{\epsilon}$ where $\bar\epsilon$ is defined in~\eqref{eq:def_bar_eps}, assume for contradiction that there exists $\epsilon_1>\Bar{\epsilon}$ such that $-p_{k}^Tg_{k}\geq \epsilon_1$. 
    Then for all $k\in\mathbb{N}$,
    \begin{equation}
    \begin{split}
        f(x_{k+1})& \leq f(x_{k})-\left[\frac{c\rho}{\alpha_0L}\left( \frac{1}{2}-\frac{\alpha_0}{2}-\left(\frac{3\alpha_0}{2}+\frac{1}{2}\right)\gamma^2 \right)\epsilon_1-\eta\right].
    \end{split}
    \end{equation}
    This shows that for each iteration there is a decrease in $f$ of at least $$\frac{c\rho}{\alpha_0L}\left( \frac{1}{2}-\frac{\alpha_0}{2}-\left(\frac{3\alpha_0}{2}+\frac{1}{2}\right)\gamma^2 \right)\epsilon_1-\eta.$$ We conclude by noting that this quantity is strictly positive as $\eps_1 >\bar\eps$ and that \\ $\frac{c\rho}{\alpha_0L}\left( \frac{1}{2}-\frac{\alpha_0}{2}-\left(\frac{3\alpha_0}{2}+\frac{1}{2}\right)\gamma^2 \right)\bar\epsilon-\eta = 0$.
    Therefore $f(x_{k})\to-\infty$ as $k\to\infty$, which is a contradiction to Assumption~\ref{assump:closed}.
    In light of Lemma~\ref{lem:uniqueness} $-p_k^Tg_k\geq 0$, and thus  $\liminf_{k\to\infty} |p_k^Tg_k|\leq \bar\epsilon$.
\end{proof}

{Finally, we remark that our analysis above can be readily extended to the self-calibrated line search, introduced in Section~\ref{sec:heuristics}.}



\section{Final Remarks}
\label{sec:finalr}

We underscore how this research extends beyond prior studies aimed at mitigating roundoff errors in optimization. 
Nonlinear optimization packages  \cite{ConnGoulToin92,GillMurrSaun05,MurtSaun83,Waec02,ZhuByrdLuNoce97} and textbooks \cite{DennSchn83,GillMurrWrig81,mybook}  devote  attention to this issue.
Nonetheless,  the strategies for handling errors are introduced as  heuristics that are seldom documented or justified. More critically, they tend to focus solely on roundoff errors,\footnote{A notable exception is Chapter 2 in  Gill, Murray and Wright \cite{GillMurrWrig81} which considers other sources of errors,  but their study is far from exhaustive.} characterized by machine precision $\epsilon_M$, which is a precisely specified quantity. 
{There is a need for a more comprehensive understanding of this topic in which stabilization techniques follow clearly specified guidelines, and where noise exhibits a more complex behavior than roundoff. This paper attempts to be a step toward that goal.}

In the future, it would be desirable to conduct similar studies, using practical applications, for more general constrained optimization problems.  We believe that the ideas presented here extend to such a wider setting.

\medskip\noindent\emph{Acknowledgments.}
The authors are grateful to Richard Byrd, Figen Oztoprak, and Stefan Wild for  valuable discussions regarding the subject matter of this paper, and to the referees for their valuable comments.
\bibliographystyle{siam}
\bibliography{references}
\newpage
\appendix
\section{Analytical Gradient of the Design Problem}  \label{sec:horn_gradient}

In the acoustic horn design problem outlined in~\eqref{prob:horn}, the gradient of the sample approximation in \eqref{precisely} can be computed analytically.
More specifically,
\begin{equation}
\tilde g_k \overset{\Delta}{=} \nabla f(b_k,\Xi_k)
= \nabla \bar{s}_k(b_k,\Xi_k)+3\nabla\sqrt{S_{k}(b_k,\Xi_k)^2},
\end{equation}
where $\bar{s}_k(b_k,\Xi_k)$ 
and $S_{k}(b_k,\Xi_k)^2$ are defined in \eqref{smalls} and~\eqref{larges}.
Furthermore, 
\be\nabla \bar{s}_k(b_k,\Xi_k) =\frac{1}{N}\sum_{\xi_i\in \Xi_k} \nabla 
s(b_k,\xi_i),\ee
and simple algebra shows that,
\begin{equation}
    \nabla \sqrt{S_k(b_k,\Xi_k)^2}=\frac{1}{N-1}\frac{\sum_{\xi_i\in \Xi_k}(s(b_k,\xi_i)-\bar{s}_k(b_k,\Xi_k))(\nabla s(b_k,\xi_i)-\nabla \bar{s}_k(b_k,\Xi_k))}{\sqrt{S_k^2(b_k,\Xi_k)}}.
\end{equation}
If $\Xi_k$ is randomly sampled, $\tilde g_k$ is an unbiased estimator for the gradient in \eqref{prob:horn} and can be computed by approximating $\nabla s(b_k,\xi_i)$ for each $\xi_i \in \Xi_k$. In the definition of $s$ from \eqref{eq:efficiency}, $\Gamma_{\text{inlet}}$ is independent of $b$. If $u$ is smooth, with $\mathds{1}$ indicating $\int_{\Gamma_{\text{inlet}}} ud\Gamma\geq0$,
\be \nabla s(b,\xi_i)= (2\mathds{1} - 1)\int_{\Gamma_{\text{inlet}}}\nabla ud\Gamma.\ee
Here $\nabla u$ can be obtained as a by-product while solving the Helmholtz equation with a finite element solver. Numerical integration over $\Gamma_{\text{inlet}}$ yields $\nabla s(b,\xi_i)$ and $\tilde g_k$. 

\section{Further Discussion on Noise Level}
\label{sec:limitations_of_noise}

In \S\ref{sec:choosing_epsA}, we  justified the rule $\epsilon_A\leftarrow\lambda \epsilon_f$, under specific assumptions on $\Delta(x)$.
However, these assumptions may not be valid in cases where the noise distribution varies significantly across different $x$ values. This limitation is evident in scenarios where canceling the mean is not possible, as discussed in \S\ref{sec:choosing_epsA}.


{While the self-calibrated strategy proposed in \S\ref{sec:heuristics}-- which can be viewed as an implicit way of estimating the local noise level-- is one possible solution, it may fail to provide sufficient safeguards
in some extreme cases (e.g. when the algorithm is highly sensitive to the choice of $\epsilon_A$). 
For those scenarios, we still need to estimate a bound on the noise $\epsilon_b$ (defined in \eqref{ebound}) or a high-probability bound. 
We discuss such estimation next.

\subsection{\boldmath Estimation of \texorpdfstring{$\epsilon_b$}{epsb} for  Stochastic Noise}
For simplicity, we will obtain $\epsilon_b$ by computing an estimate of $\sup \|\Delta(x)\|$ at a representative $x$.
A global estimate can then be derived  e.g. by \eqref{eq:epsf_avg}.

One can establish consistent estimators of the noise bound if we can compute an estimate on the true objective value.
Let us generate $m$ i.i.d. samples \\  $\{\tilde f_1(x), \tilde f_2(x),\ldots,\tilde f_m(x)\}$ and let us compute an accurate estimate of the true objective $f(x)$, denoted as $\widehat f(x)$. 
Then the samples of noise in the function are given by
    \begin{equation}
        \delta_j(x):=\tilde f_j(x)-\widehat f(x),\quad j=1,2,\cdots,m.
    \end{equation}  
A concrete example arises in stochastic optimization where the true objective is $f(x):=\mathbb{E}(F(x,\xi))$.
The $j$th sample of the noisy objective is defined as  $\tilde f_j(x)=\frac{1}{N}\sum_{i=1}^N F(x,\xi_{j_i})$ for an i.i.d. batch $\{\xi_{j_1},\xi_{j_2},\cdots,\xi_{j_N}\}$ of size $N$.
An accurate estimator $\widehat f(x)$ of $f(x)$ can then be defined as $\widehat f(x)=\frac{1}{M}\sum_{i=1}^M F(x,\xi_{i})$ for another batch of i.i.d. samples $\{\xi_{i}\}_{i=1}^M$, where $M\gg N$ is sufficiently large.

We provide the following three estimators that can be used in practice, where the first two require the access to $\widehat f(x)$ and the third one does not: 
\medskip

    1) \emph{Empirical Chebyshev bound \cite{saw1984chebyshev}:}
    \begin{equation}
        \hat{\epsilon}_b^1:=\overline{\delta(x)}+\lambda \sqrt{\frac{1}{m-1}\sum_{j=1}^m(\delta_j(x)-\overline{\delta(x)})^2}
    \end{equation}
    for some integer $\lambda$ large enough, where $\overline{\delta(x)}= [ \delta_1(x)+\cdots+\delta_m(x)]/m$ .
    
    2) \emph{Maximum of $|\delta_j(x)|$:}
    \begin{equation}
        \hat{\epsilon}_b^2:=\max_{j=1,\cdots,m}\{ |\delta_j(x)| \}.
    \end{equation}

    3) \emph{Range of noisy objectives:}
    \begin{equation}
    \label{eq:range_estimator}
        \hat{\epsilon}_b^3:=\max_{j=1,\cdots,m} \tilde f_j(x)-\min_{j=1,\cdots,m} \tilde f_j(x).
    \end{equation}

$\hat{\epsilon}_b^1$ is a high-probability bound of $\|\Delta(x)\|$, assuming that the noise has a finite variance but not necessarily bounded.
$\hat{\epsilon}_b^2$ is a consistent estimator of $\epsilon_b$ if $\sup\|\Delta(x)\|<\infty$.
$\hat{\epsilon}_b^3$ can be a biased (and depending on the estimated quantity, potentially inconsistent) estimator if the noise does not have mean zero, yet it can be easily computed without $\widehat{f}(x)$.
In practice, $\hat{\epsilon}_b^3$ is an attractive candidate when $\widehat{f}(x)$ is expensive or not accessible, or when the noise level estimate is not required to be accurate, as in the acoustic horn design.


\subsection{\boldmath \texorpdfstring{Estimating $\epsilon_b$}{epsb} for Computational Noise}
Due to the deterministic nature of computational noise, the first two estimators discussed above cannot be employed.
As an alternative, we can modify the range estimator~\eqref{eq:range_estimator} following a similar approach as \texttt{ECNoise}.
At a selected point $x$, one can collect noisy objectives in a  small neighborhood of $x$, and then compute the range as an estimate of $\eps_b$.
Similar to the argument for stochastic noise, if the distribution does not vary significantly, using $\texttt{ECNoise}$ is usually effective; see \cite{more2011estimating}.

\section{Sample Selection and Consistency}
\label{sec:sample_consistency}
In many stochastic optimization problems, such as the acoustic horn design described in \S\ref{sec:case}, the noisy evaluations $\tilde f(x_k)$, e.g.~\eqref{precisely}, depend on a particular sample batch $\Xi_k$. 
In certain cases, the selection of $\Xi_k$ is entirely under the control of the user. One can thus fix $\Xi_k$ during the course of an iteration of the optimization algorithm,  a case we  refer to as ``sample consistency''. In such a setting the effect of noise on function comparisons and differences is more benign.


Reusing samples is, however, not always possible. In that case, the algorithm will operate in the ``sample inconsistent'' regime, which is the most general and challenging for optimization methods and holds particular interest in this paper. 

Let us summarize these two cases for the key components of our algorithm.

\paragraph{Relaxed line search} For backtrack numbers $\ell=1,2, \cdots$, we denote the sample used in the evaluation of $\tilde f(x_k+\beta_k^\ell \tpk)$ by $\Xi_k^\ell$.
In the sample inconsistent case, the $\Xi_k^\ell$ are different from each other and a relaxation $\epsilon_A$ is employed.
On the other hand, if sample consistency is ensured, we can set $\epsilon_A\leftarrow 0$ since no errors are involved in the comparison with a fixed $\Xi_k$
\footnote{Note that although the comparison is robust, $\tilde f(x_k)$ is still a noisy estimate and a careful choice of $\alpha_0$ can be useful when noise is large; see \S\ref{sec:sample_consistent_numerical_results}.}.

\paragraph{Finite differences} Given the estimated noise level $\epsilon_f$, the finite difference estimator is
\begin{equation}
\label{eq:FDHorn}[\tilde g_k]_i:=\frac{\tilde f(b_k+he_i,\Xi_k^2)-\tilde f(b_k,\Xi_k^1)}{h}\qquad i=1,\cdots,n,
\end{equation}
where $\Xi_k^1$ and $\Xi_k^2$ are two batches.
Sample inconsistency allows $\Xi_k^1\neq \Xi_k^2$, and $h$ needs to be chosen according to the noise level as seen in \eqref{eq:FDOptimalMSE}.
With sample consistency, $\Xi_k^1=\Xi_k^2$, formula \eqref{eq:FDHorn} gives a fairly accurate gradient approximation of the corresponding sample average approximation of the objective, and thus $h$ is set as the unit roundoff $\epsilon_M$.


\subsection{Numerical Results with Sample Consistency}
\label{sec:sample_consistent_numerical_results}

We study the performance of algorithm {\tt GP-LS} when fixing the sample  during line search and gradient estimation.
In Figure~\ref{fig:samp_consis}, we plot the performance of {\tt GP-LS} with $\epsilon_A=0$, and  for  $N = 10, 50, 100$.
For each value of $N$, we adjust $\alpha_0$ ($0.1,0.25,1$ respectively) to cope with the fact that the sample average approximations of the objective function become increasingly inaccurate as $N$ decreases.
The finite difference interval $h$ is chosen to be $10^{-6}$ for all cases.

\begin{figure}[!htbp]
    \centering
    \includegraphics[width = 0.49\textwidth]{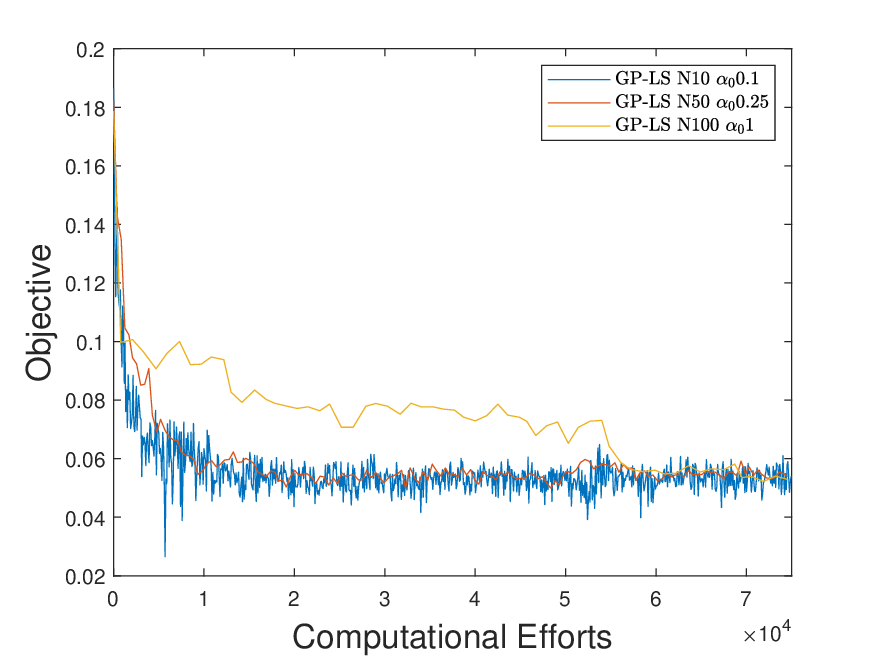}
    \includegraphics[width = 0.49\textwidth]{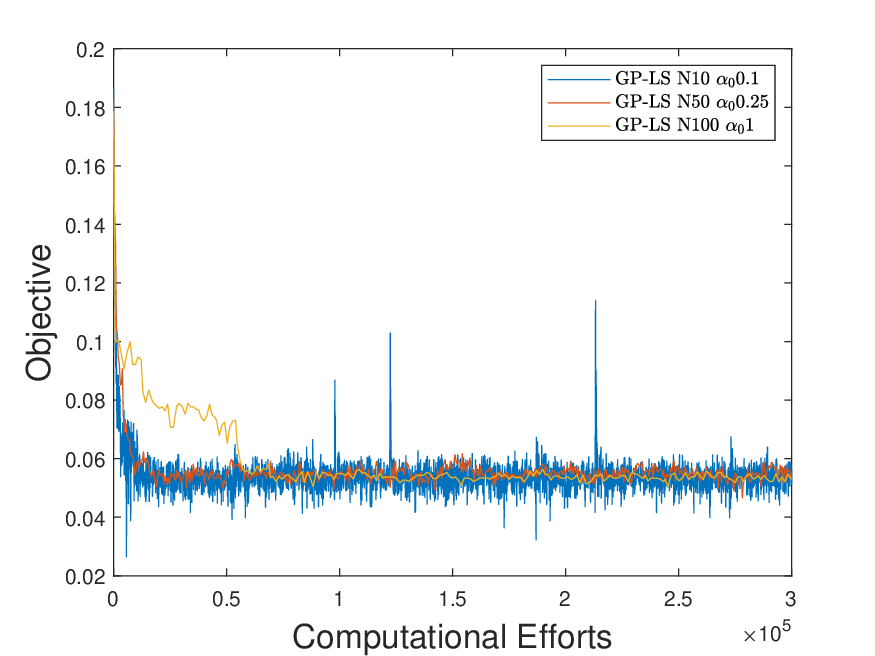}
    \caption{Comparison of different sample sizes when using a sample consistent version of Algorithm {\tt GP-LS} using sample sizes $10,50,100$, and different $\alpha_0$ respectively;  Left: Objective function value vs. computational effort (up to $75,000$). Right: Objective function value vs. computational effort (up to $3\times 10^5$).}
    \label{fig:samp_consis}
\end{figure}

We observe in Figure~\ref{fig:samp_consis} that all three plots exhibit nice convergence behavior. With smaller sample sizes the iterates approach the solution more quickly, although they may give rise to spikes as the iteration continues.
We conclude that, when feasible, sample consistency results in robust and efficient performance, if an appropriate value of the sample size $N$ is first determined after experimentation.

\section{Supplementary Proof}
\label{sec:supp_proof}
\begin{lemma}
\label{lem:StatExistApdx}
     $x^*\in \Omega$ is a first-order stationary point of problem \eqref{opt_problem} if and only if $p(x^*)=0$.
\end{lemma}

\begin{proof}
    Prop. 6.1.1 (b) in \cite{bertsekas2015convex} shows $(\Leftarrow)$ of Lemma~\ref{lem:StatExist}.
    
    To see $(\Rightarrow)$, since $x^*$ is a stationary point and by definition, $ g(x^*)^T(x-x^*)\geq 0$ for all $x\in\Omega$.
    Take $x=P_\Omega[x^*-\alpha_0 g(x^*)]$, then
    \begin{equation}
    \label{ineq:stat_point_prod_proj}
         g(x^*)^T(P_\Omega[x^*-\alpha_0 g(x^*)]-x^*)=p(x^*)^Tg(x^*)\geq 0.
    \end{equation}
    Note that by letting $x=x^*-\alpha_0 g(x^*)$, $z=P_{\Omega}[x^*-\alpha_0 g(x^*)]$ and $y=x^*$ in Lemma~\ref{lem:uniqueness}, one has
    \begin{equation*}
    \begin{split}
        & \left( x^*-\alpha_0 g(x^*) -P_{\Omega}[x^*-\alpha_0 g(x^*)] \right)^T\left( x^*-P_{\Omega}[x^*-\alpha_0 g(x^*)] \right)\leq 0 \\
        \Longrightarrow  & \|p(x^*)\|^2=\|x^*-P_{\Omega}[x^*-\alpha_0 g(x^*)]\|^2\leq -\alpha_0 p(x^*)^Tg(x^*)\leq 0
    \end{split}
    \end{equation*}
    where the final inequality follows from~\eqref{ineq:stat_point_prod_proj}.
    This implies that $p(x^*)=0$.
\end{proof}

\end{document}